\documentclass[12pt]{article}
\usepackage{latexsym}
\usepackage{a4wide}
\usepackage{theorem}
\usepackage{amsfonts}
\usepackage{enumerate}
\usepackage{amsmath}
\usepackage{amscd}
\usepackage{times}

\theoremstyle{change}
\theorembodyfont{\sl}
\newtheorem{thm}[subsection]{Theorem.}
\newtheorem{prop}[subsection]{Proposition.}

\newtheorem{cor}[subsection]{Corollary.}
\newtheorem{conj}[subsection]{Conjecture.}
\newtheorem{sublem}[subsection]{Lemma.}
\newtheorem{subthm}[subsection]{Theorem.}
\newtheorem{subprop}[subsection]{Proposition.}

\theorembodyfont{\rmfamily}
\newtheorem{defi}[subsection]{Definition.}
\newtheorem{rem}[subsection]{Remark.}

\makeatletter
\newenvironment{eqn}{\refstepcounter{subsection}
$$}{\leqno{\rm(\thesubsection)}$$\global\@ignoretrue}
\makeatother  







\newenvironment{prf}[1]{\trivlist
\item[\hskip\labelsep{\it
#1.\hspace*{.3em}}]}{~\hspace{\fill}~$\square$\endtrivlist}

\newenvironment{proof}{\begin{prf}{\bf Proof}}{\end{prf}}

\newcommand{\square}{\Box}
\newcommand{\lto}{\longrightarrow}

\newcommand{\into}{\hookrightarrow}
\newcommand{\ZZ}{{\mathbb Z}}
\newcommand{\FF}{{\mathbb F}}
\newcommand{\QQ}{{\mathbb Q}}
\newcommand{\PP}{{\mathbb P}}
\newcommand{\RR}{{\mathbb R}}
\newcommand{\CC}{{\mathbb C}}
\newcommand{\Zhat}{\hat{\ZZ}}
\renewcommand{\AA}{{\mathbb A}}
\newcommand{\AAf}[1][]{\AA_{{#1}{\rm f}}}
\renewcommand{\SS}{{\mathbb S}}
\newcommand{\HH}{{\mathbb H}}
\newcommand{\Spec}{{\rm Spec}}
\newcommand{\Pic}{{\rm Pic}}
\newcommand{\GG}{{\mathbb G}}
\newcommand{\Gm}[1][]{{\GG_{{\rm m}#1}}}
\newcommand{\eps}{\varepsilon}
\newcommand{\id}{{\rm id}}
\newcommand{\GL}{{\rm GL}}
\newcommand{\PGL}{{\rm PGL}}
\newcommand{\SL}{{\rm SL}}
\newcommand{\GSp}{{\rm GSp}}
\newcommand{\trace}{{\rm tr}}
\newcommand{\Aut}{{\rm Aut}}
\newcommand{\ol}{\overline}
\newcommand{\calO}{{\cal O}}
\newcommand{\calL}{{\cal L}}
\newcommand{\Hom}{{\rm Hom}}
\newcommand{\End}{{\rm End}}
\newcommand{\Ext}{{\rm Ext}}
\newcommand{\Sym}{{\rm Sym}}
\newcommand{\ld}{\langle}
\newcommand{\rd}{\rangle}
\newcommand{\Sh}{{\rm Sh}}
\newcommand{\Vol}{{\rm Vol}}
\newcommand{\rH}{{\rm H}}
\newcommand{\rM}{{\rm M}}
\newcommand{\Res}{{\rm Res}}
\newcommand{\MT}{{\rm MT}}
\newcommand{\der}{{\rm der}}
\newcommand{\swedge}{{\scriptstyle{\wedge}}}

\newcommand{\Gal}{{\rm Gal}}
\newcommand{\Qbar}{{\ol{\QQ}}}
\newcommand{\ad}{{\rm ad}}
\newcommand{\intaut}{{\rm int}}
\newcommand{\AM}{{\rm AM}}
\newcommand{\sm}{{\rm sm}}
\newcommand{\nor}{{\rm nor}}
\newcommand{\discr}{{\rm discr}}
\newcommand{\ab}{{\rm ab}}
\newcommand{\Norm}{{\rm Norm}}
\newcommand{\Reg}{{\rm Reg}}
\newcommand{\tors}{{\rm tors}}
\newcommand{\et}{{\rm et}}
\newcommand{\Li}{{\rm Li}}
\newcommand{\Frob}{{\rm Frob}}

\begin{document}
\title{On the Andr\'e-Oort conjecture for Hilbert modular surfaces.
\thanks{A preprint version of this article together with an appendix
that could be considered as the author's scratch paper while working
on this subject can be downloaded from the author's home page. This
appendix contains details that the author does not find interesting
enough to publish, but that may be helpful for readers who got lost.}}
\author{Bas Edixhoven\footnote{partially supported by the Institut
Universitaire de France, and by the European TMR Network Contract ERB
FMRX 960006 ``arithmetic algebraic geometry''.}}  
\date{April 5, 2000}

\maketitle

\section{Introduction.}\label{sec1}
In order to state the conjecture mentioned in the title, we need to
recall some terminology and results on Shimura varieties; as a general
reference for these, we use~\cite[Sections~1--2]{Moonen2}.  So let
$\SS:=\Res_{\CC/\RR}\GG_{{\rm m},\CC}$ be the algebraic group over
$\RR$ obtained by restriction of scalars from $\CC$ to $\RR$ of the
multiplicative group. For $V$ an $\RR$-vector space, it is then
equivalent to give an $\RR$-Hodge structure or an action by $\SS$ on
it. A Shimura datum is a pair $(G,X)$, with $G$ a connected reductive
affine algebraic group over $\QQ$, and $X$ a $G(\RR)$-conjugacy class
in the set of morphisms of algebraic groups $\Hom(\SS,G_\RR)$,
satisfying the three conditions of \cite[Def.~1.4]{Moonen2} (i.e., the
usual conditions (2.1.1--3) of~\cite{Deligne1}). These conditions
imply that $X$ has a natural complex structure (in fact, the connected
components are hermitian symmetric domains), such that every
representation of $G$ on a $\QQ$-vector space defines a polarizable
variation of Hodge structure on~$X$. For $(G,X)$ a Shimura datum, and
$K$ a compact open subgroup of $G(\AAf)$, we let $\Sh_K(G,X)(\CC)$
denote the complex analytic variety $G(\QQ)\backslash(X\times
G(\AAf)/K)$, which has a natural structure of quasi-projective complex
algebraic variety, denoted $\Sh_K(G,X)_\CC$; the projective limit
$\Sh(G,X)_\CC$ over all $K$ of the $\Sh_K(G,X)_\CC$ is a scheme on
which $G(\AAf)$ acts continuously. (The action being continuous means
that the scheme has a cover by open affines $U_i=\Spec(A_i)$ such that
each $U_i$ is stabilized by some open subgroup $K_i$ of $G(\AAf)$ and
each $f$ in $A_i$ has open stabilizer in~$K_i$.)  A morphism of
Shimura data from $(G_1,X_1)$ to $(G_2,X_2)$ is a morphism $f\colon
G_1\to G_2$ that maps $X_1$ to $X_2$; for $K_1$ and $K_2$ compact open
subgroups of $G_1(\AAf)$ and $G_2(\AAf)$ with $f(K_1)$ contained in
$K_2$, such an $f$ induces a morphism $\Sh(f)$ from
$\Sh_{K_1}(G_1,X_1)_\CC$ to $\Sh_{K_2}(G_2,X_2)_\CC$.
\begin{defi}\label{def1.1}
Let $(G,X)$ be a Shimura datum, $K$ an open compact subgroup of
$G(\AAf)$, and $Z$ an irreducible closed subvariety of
$\Sh_K(G,X)_\CC$. Then $Z$ is a {\sl subvariety of Hodge type} if
there is a Shimura datum $(G',X')$, a morphism of Shimura data
$f\colon(G',X')\to(G,X)$, and an element $g$ of $G(\AAf)$ such that
$Z$ is an irreducible component of the image of the map:
$$
\Sh(G',X')_\CC \overset{\Sh(f)}{\lto} \Sh(G,X)_\CC \overset{{\cdot}g}{\lto} 
\Sh(G,X)_\CC 
\lto \Sh_K(G,X)_\CC. 
$$
\end{defi}
This definition is equivalent to \cite[6.2]{Moonen2}, which uses only
closed immersions $f\colon G'\to G$.  In \cite[Prop.~2.8]{Moonen3} it
is explained that the subvarieties of Hodge type are precisely the
loci where certain given classes in certain variations of Hodge
structures (obtained from representations of $G$) are Hodge classes;
hence the terminology.
\begin{defi}
Let $(G,X)$ be a Shimura datum. For $h$ in $X$ we let $\MT(h)$ be the
Mumford-Tate group of $h$, i.e., the smallest algebraic subgroup $H$
of $G$ such that $h$ factors through~$H_\RR$. A point $h$ in $X$ is
called {\sl special} if $\MT(h)$ is commutative (in which case it is a
torus).  For $K$ a compact open subgroup of $G(\AAf)$, a point in
$\Sh_K(G,X)_\CC$ is special if its preimages in $\Sh(G,X)_\CC$ are of
the form $(h,g)$ with $h$ in $X$ special. Equivalently, the special
points in $\Sh_K(G,X)_\CC$ are the zero dimensional subvarieties of
Hodge type.
\end{defi}
\begin{conj}[Andr\'e-Oort]\label{conj1.3}
Let $(G,X)$ be a Shimura datum. Let $K$ be a compact open subgroup
of~$G(\AAf)$ and let $S$ be a set of special points in
$\Sh_K(G,X)(\CC)$. Then every irreducible component of the Zariski
closure of $S$ in $\Sh_K(G,X)_\CC$ is a subvariety of Hodge type.
\end{conj}
Some remarks are in order at this point. Andr\'e stated this
conjecture as a problem for curves containing infinitely many special
points in general Shimura varieties in \cite[X.4]{Andre1}.
Independently, Oort raised the question for general subvarieties of
the moduli spaces of principally polarized abelian varieties
in~\cite{Oort1}.  In \cite[X.4]{Andre1} Andr\'e mentions the
similarity with the Manin-Mumford conjecture (proved by Raynaud,
see~\cite{Abbes1}), and \cite{Andre2} contains a version generalizing
both the conjecture above and the Manin-Mumford conjecture; see also
\cite[6.7.2]{Moonen2}.

Let us now discuss the results on the conjecture above that have been
obtained until now. All of them deal with moduli spaces of abelian
varieties. Moonen proved in his thesis (see \cite[\S5]{Moonen4}, in
particular the equivalence between Conjectures~5.1 and 5.3, and
\cite[IV]{Moonen1}) that the conjecture is true for sets $S$ for which
there exists a prime number $p$ at which all $s$ in $S$ have an
ordinary reduction of which they are the canonical lift. Needless to
say, his methods use reduction modulo a prime number~$p$.  This gives
a quite general result, but it has the disadvantage that one neglects
most of the Galois action on the special points, and that one has to
work with one Frobenius element simultaneously for all $s$ in~$S$.

In \cite{Edixhoven1}, the conjecture was proved for the moduli space
of pairs of elliptic curves, assuming the generalized Riemann
hypothesis (GRH) for imaginary quadratic fields. In a few words, the
proof exploits the Galois action on CM-points and considers
intersections of the subvarieties in question with images of them
under suitable Hecke operators. In this approach, we work with a
different Frobenius element for each $s$ in~$S$; GRH comes in via the
existence of small primes with suitable properties.  The same case of
Conjecture~\ref{conj1.3} was proved unconditionally by Andr\'e
in~\cite{Andre3}. He uses the Galois action on the CM-points, and a
Diophantine approximation result of Masser on the $j$-function.

More recently, Yafaev has generalized the result in \cite{Edixhoven1}
to the case of products of two Shimura curves that are associated to
quaternion algebras over $\QQ$, see~\cite{Yafaev1}, and B.~Belhaj
Dahman, a student of Andr\'e, is working on the families of jacobians
of the curves
$$
y^n=x(x-1)(x-\lambda). 
$$
The question about these families of jacobians is whether or not the 
various isogeny factors coming from the decomposition for the action of 
$\mu_n(\CC)$ are simultaneously of CM type for infinitely many complex 
numbers~$\lambda$. 

Recently, Clozel and Ullmo have proved (\cite{ClozelUllmo1}), for $G$
among $\GSp_{2n}$ and $\GL_n$, that sets of the form $T_px$, with $x$
in $G(\QQ)\backslash G(\AA)/K$ and $T_p$ certain Hecke operators with
$p$ tending to infinity, are equidistributed.  The idea behind this is
that one would like to imitate and apply the equidistribution results
for Galois orbits of strict sequences of points of small height in
abelian varieties as in~\cite{Abbes1}. A sequence of closed points of
an algebraic variety is called strict if every proper closed subset
contains only finitely many elements of the sequence. Phrased in this
terminology, the Andr\'e-Oort conjecture says that a sequence of
special points is strict if no proper subvariety of Hodge type
contains an infinite subsequence.  Of course, to prove the
Andr\'e-Oort conjecture in this way, one has to replace $T_px$ by the
Galois orbit of $x$, which seems to be a hard problem, and moreover,
one has to deal with the fact that the heights of CM points tend to
infinity and not to zero.

In this article, we prove the Conjecture~\ref{conj1.3}, assuming GRH,
for Hilbert modular surfaces. The method of proof is basically the
same as in~\cite{Edixhoven1}, but now we do use more advanced
techniques. The two main results of the article are described in
Section~\ref{sec3}. The reason for which we state and prove
Theorem~\ref{thm3.2} is that it has an interesting application to
transcendence of special values of certain hypergeometric functions
via work of Wolfart, Cohen and W\"ustholz, see~\cite{CohenWustholz1},
{\sl without} having to assume GRH.

Let us briefly describe the contents of this
article. Section~\ref{sec3} introduces the Hilbert modular surfaces
that we work with in terms of a Shimura datum, gives their
interpretation as moduli spaces of abelian surfaces with
multiplications by the ring of integers of a real quadratic field $K$,
and states the main results.

Section~\ref{sec3.3}, which is not so essential, discusses the
difference between working with abelian surfaces with or without a
given polarization. In group theoretical terms, the choice is between
working with $\GL_2(K)$ or its subgroup $\GL_2(K)'$ consisting of the
elements of $\GL_2(K)$ whose determinant is in~$\QQ^*$. The reason for
considering both cases is that with a polarization (and a suitable
level structure), the variation of Hodge structure provided by the
lattices of the abelian surfaces comes from a representation of the
group in the Shimura datum, which is not true without given
polarizations. We need variations of Hodge structure in
Section~\ref{sec3.4}. On the other hand, the size of Galois orbits of
special points, studied in Section~\ref{sec3.6}, is simpler to
understand in terms of class groups when working without
polarizations. We could have chosen to work throughout the article
with $\GL_2(K)'$, but we think that it is instructive to see the
consequences of such a choice in the relatively easy case of Hilbert
modular surfaces, before trying to treat general Shimura varieties
completely in group theoretical terms.

In Section~\ref{sec3.4} we recall an important result of Andr\'e,
relating the generic Mumford-Tate group of a variation of Hodge
structure to its algebraic monodromy group (i.e., the Zariski closure
of the image of monodromy). We use it to prove that for a curve in a
Hilbert modular surface that is not of Hodge type and that does
contain a special point, the connected algebraic monodromy group is
maximal, i.e.,~$\SL_{2,K}$.

Section~\ref{sec3.5} introduces the Hecke correspondence $T_p$
associated to a prime number~$p$. We use a very powerful result of
Nori in order to prove that for $C$ a curve with maximal algebraic
monodromy group, $T_pC$ is irreducible if $C$ is large enough.

The main result of Section~\ref{sec3.6} says that the size of the
Galois orbit of a special point $x$ grows at least as a positive power
of the discriminant $\discr(R_x)$ of the ring of endomorphisms
(commuting with the real multiplications) of the corresponding abelian
variety. This section is quite long, and contains some messy
computations, depending on the structure of the Galois group of the
normal closure of the CM field in question. The problem is that one
has to give a lower bound for the image under the reflex type norm of
one class group in another.

Section~\ref{sec3.7} gives an upper bound for the number of points in
intersections of the form $Z_1\cap T_gZ_2$, with $Z_1$ and $Z_2$ fixed
subvarieties of a general Shimura variety, and with $T_g$ a varying
Hecke correspondence.

Finally, Section~\ref{sec3.8} combines all these preliminary results
as follows. One supposes that $C$ is a curve in a Hilbert modular
surface $S$, containing infinitely many special points, and not of
Hodge type. If $p$ is large enough (depending only on $C$), then
$T_pC$ is irreducible by Section~\ref{sec3.5}. Since the $T_p$-orbits
in $S$ are dense, one cannot have $C=T_pC$. Hence the intersections
$C\cap T_pC$ are finite, and hence bounded above
(Section~\ref{sec3.7}) by a constant times~$p^2$. Let now $x$ be a
special point on~$C$. If $p$ is a prime that is split in $R_x$, then
$C\cap T_pC$ contains the Galois orbit of~$x$, hence $|C\cap T_pC|$
grows at least as a positive power of~$|\discr(R_x)|$. But this lower
bound for primes that are split in $R_x$ contradicts the conditional
effective Chebotarev theorem (this is where GRH comes in). Hence,
assuming GRH, one has proved that if $C$ does contain infinitely many
special points, then $C$ is of Hodge type. The reason that one can
prove Thm.~\ref{thm3.2} unconditionally is that in that case the CM
field $\QQ\otimes R_x$ is independent of $x$, and hence Chebotarev's
theorem itself is sufficient.

In April 1999, we have proved Conjecture~\ref{conj1.3}, assuming GRH,
for arbitrary products of modular curves, extending the methods
of~\cite{Edixhoven1}. A detailed proof, which is quite elementary,
will be written up in the near future. One can hope that combining the
techniques used for these last two results will make it possible to
treat more general higher dimensional cases of
Conjecture~\ref{conj1.3}. Of course, eventually everything should be
expressed in terms of ``$(G,X)$-language''.  In fact, in this article
we could already have worked without mentioning abelian varieties.

Before we really start, let us first mention two obvious general
principles.  The first is that level structures don't matter in
Conjecture~\ref{conj1.3}: for $(G,X)$ a Shimura datum, $K$ and $K'$
open compact in $G(\AAf)$ with $K\subset K'$, an irreducible
subvariety $Z$ of $\Sh_K(G,X)_\CC$ is of Hodge type if and only if its
image in $\Sh_{K'}(G,X)_\CC$ is. The second principle says that the
irreducible components of intersections of subvarieties of Hodge type
are again of Hodge type (this is clear from the interpretation of
subvarieties of Hodge type given right after Definition~\ref{def1.1}).

\section{The main results.}\label{sec3}
Let $K$ be a real quadratic extension of $\QQ$, let $O_K$ be its ring
of integers, and let $G$ be the $\ZZ$-group scheme
$\Res_{O_K/\ZZ}(\GL_{2,O_K})$. After numbering the two embeddings of
$K$ in $\RR$, we have $\RR\otimes K=\RR^2$, and hence
$G(\RR)=\GL_2(\RR)^2$. We will study the Shimura variety:
$$
S(\CC):=G(\QQ)\backslash (X\times G(\AAf)/G(\Zhat)), 
$$
where $X=(\HH^{\pm})^2$, and where $\HH^\pm$ is the usual
$\GL_2(\RR)$-conjugacy class of morphisms from $\SS$ to $\GL_{2,\RR}$,
i.e., the class of $a+bi\mapsto\bigl(\begin{smallmatrix}a & -b\\b &
a\end{smallmatrix}\bigr)$.  The surface $S_\CC$, called a Hilbert
modular surface, is the coarse moduli space for pairs $(A,\alpha)$
with $A$ an abelian surface and $\alpha$ a morphism from $O_K$ to
$\End(A)$ (see~\cite[Ch.~X]{vanderGeer1}, and the end of
Section~\ref{sec3.3} for the moduli interpretation for a closely
related Shimura datum). This implies that the reflex field of $(G,X)$
is $\QQ$ and that the canonical model $S_\QQ$ (see
\cite[Section~2]{Moonen2} for this notion) is simply the coarse moduli
space for pairs $(A/S/\QQ,\alpha)$ with $S$ a $\QQ$-scheme, $A/S$ an
abelian scheme of relative dimension two, and $\alpha$ a morphism from
$O_K$ to~$\End_S(A)$.  The set of geometrically connected components
of $S_\QQ$ is $K^*\backslash\AA_K^*/(\RR\otimes K)^{*,+}O_K^{\swedge
*}=\Pic(O_K)^+$, the group of isomorphism classes of invertible
$O_K$-modules with orientations at the two infinite places, and has
trivial action by $\Gal(\Qbar/\QQ)$ (\cite[Ch.~I,
Cor.~7.3]{vanderGeer1}).  The main objective of this article is to
prove the following two theorems.
\begin{thm}\label{thm3.1}
Assume GRH. Let $C\subset S_\CC$ be an irreducible closed curve 
containing infinitely many CM points. Then $C$ is of Hodge type. 
\end{thm}
\begin{thm}\label{thm3.2}
Let $C\subset S_\CC$ be an irreducible closed curve containing
infinitely many CM points corresponding to abelian varieties that lie
in one isogeny class (the isogenies are not required to be compatible
with the multiplications by~$O_K$). Then $C$ is of Hodge type.
\end{thm}
Let us note immediately that these theorems apply in fact to all
Hilbert modular surfaces, because the Andr\'e-Oort conjecture is
insensitive to level structure. Before proving the theorems we need to
discuss some of the tools we will use in it.

\section{Choosing a suitable Shimura variety.}\label{sec3.3}
For a variation of Hodge structure on a complex variety, one has the
notions of generic Mumford-Tate group and that of monodromy. A
relation between these two notions will be very useful for us. In
order to get a suitable variation of Hodge structure on $S(\CC)$ as
above, there is a little complication, and at least two options to get
around it. The problem is that the tautological representation of
$G_\QQ$ on the $\QQ$-vector space $K^2$ does not induce a variation of
Hodge structure on $\Sh_H(G_\QQ,X)(\CC)$, even if $H$ is an arbitrary
small open subgroup of $G(\AAf)$; just consider the action of $O_K^*$
in $G(\QQ)$ on $X\times G(\AAf)/H$ (see \cite[Section~2.3]{Moonen3}
for a general statement).

The first possible way out is to use an other representation, and have
the monodromy take place in the image of $G$ under this
representation.  For example, one can take the representation
$\Sym^2(\rho_0)\otimes\det(\rho_0)^{-1}$, with $\rho_0$ the
tautological representation on~$O_K^2$. This representation induces a
faithful representation $\rho$ of the quotient
$G^\ad=\Res_{O_K/\ZZ}\PGL_{2,O_K}$.  The morphism $G\to G^\ad$ induces
an isomorphism from $X$ to a conjugacy class $X^\ad$ in
$\Hom_\RR(\SS,G^\ad_\RR)$, and gives a morphism of Shimura data from
$(G,X)$ to $(G^\ad,X^\ad)$. Let $S^\ad(\CC)$ be the Shimura variety
$\Sh_{G^\ad(\Zhat)}(G^\ad,X^\ad)(\CC)$. The natural morphism from
$S_\CC$ to $S^\ad_\CC$ is finite and surjective (this follows directly
from the definition); one can show that it is the quotient for a
faithful action of~$\Pic(O_K)$, but that will not be used.
Conjecture~\ref{conj1.3} is then true for $S_\CC$ if and only if it is
true for $S^\ad_\CC$, and $\rho$ induces a variation of Hodge
structure on $\Sh_H(G^\ad,X^\ad)(\CC)$ for suitable~$H$.  The
disadvantage of working with $S^\ad_\CC$ is that it does not seem to
have an interpretation as a moduli space of abelian varieties; this is
not a real problem, but we prefer to work with Shimura varieties that
are as simple as possible.

Another way out is to replace the group $G$ by its subgroup $G'$ given
by the following Cartesian diagram:
$$
\begin{array}{ccc}
G' & \into & G \\
\downarrow & \square & \downarrow\scriptstyle{\det} \\
\Gm[,\ZZ] & \into & \Res_{O_K/\ZZ}\Gm[,O_K] 
\end{array}
$$
Loosely speaking, $G'$ is the subgroup of $G$ consisting of those
elements whose determinant is in~$\QQ$. As the morphism $\det$ in the
diagram above is smooth, $G'$ is smooth over $\Gm[,\ZZ]$, hence
over~$\ZZ$.  It follows that $G'$ is the scheme-theoretic closure in
$G$ of its generic fibre. We note that $G'(\RR)$ is the subgroup of
$(x,y)$ in $\GL_2(\RR)^2$ with $\det(x)=\det(y)$.  All $h$ in $X$
factor through $G'_\RR$, but $X$ consists of two $G'(\RR)$-conjugacy
classes. The conjugacy class $X'$ we work with is the disjoint union
of $(\HH^+)^2$ and~$(\HH^-)^2$. This gives a morphism of Shimura data
from $(G',X')$ to $(G,X)$, and a morphism of Shimura varieties
$S'_\CC\to S_\CC$ with $S'_\CC=\Sh_{G'(\Zhat)}(G',X')_\CC$. One can
prove that the Shimura variety $S'_\CC$ is connected, and that the
morphism to its image in $S_\CC$ is the quotient by a faithful action
of the finite group $O_K^{*,+}/O_K^{*,2}$, i.e., by the group of
totally positive global units modulo squares of global units. We will
only use that the morphism $S'_\CC\to S_\CC$ is finite and that its
image is open and closed; these two facts follow directly from the
definitions. It follows that Conjecture~\ref{conj1.3} is true for
$S_\CC$ if and only if it is for~$S'_\CC$, and similarly for the two
theorems above that we want to prove. Moreover, the tautological
representation of $G'$ does induce a variation of Hodge structure on
$\Sh_H(G',X')(\CC)$ for $H$ sufficiently small.

The option we choose is the last. The variety $S'_\CC$ is the (coarse)
moduli space for triplets $(A,\alpha,\lambda)$ where: 
\begin{eqn}\label{eqn3.1}
\left\{ 
\begin{array}{l}
\text{$A$ is a complex abelian surface,}\\
\text{$\alpha\colon O_K\to \End(A)$ a morphism of rings,}\\
\text{and $\lambda\colon A\to A^*$ a principal $O_K$-polarization,}
\end{array}
\right.
\end{eqn} a notion that we will now explain. First of
all, $A^*$ is the dual of $A$ in the category of abelian varieties
with $O_K$-action: $A^*:=\Ext^1(A,O_K\otimes\Gm)$. One verifies that
$A^*=\delta\otimes_{O_K}A^t$, where $\delta$ is the different of the
extension $\ZZ\to O_K$, and where $A^t=\Ext^1(A,\Gm)$, the dual of $A$
in the usual sense.  The inclusion $\delta\subset O_K$ induces a
morphism $A^*\to A^t$, which is an isogeny. A principal
$O_K$-polarization is then an isomorphism $\lambda\colon A\to A^*$
such that the induced morphism from $A$ to $A^t$ is a
polarization. Interpreted in Hodge-theoretical terms, a triplet
$(A,\alpha,\lambda)$ as in (\ref{eqn3.1}) corresponds to a triplet
$(V,h,\psi)$ with $V$ a locally free $O_K$-module of rank two,
$h\colon\SS\to(\GL_\ZZ(V))_\RR$ a Hodge structure of type
$(-1,0),(0,-1)$, and $\psi\colon V\times V\to O_K$ a perfect
antisymmetric $O_K$-bilinear form such that $\trace\circ\psi\colon
V\times V\to\ZZ$ is a polarization.  Note that for such a triplet
$(V,h,\psi)$, the pair $(V,\psi)$ is isomorphic to the standard pair
$(O_K\oplus O_K,(\begin{smallmatrix}0&1\\-1&0\end{smallmatrix}))$. In
order to prove that the set of isomorphism classes of
$(A,\alpha,\lambda)$ as in (\ref{eqn3.1}) is $S'(\CC)$ one uses the
following two facts: 1:~$G'(\AAf)/G'(\Zhat)$ is the set of
$O_K$-lattices in $K^2$ on which
$\psi=(\begin{smallmatrix}0&1\\-1&0\end{smallmatrix})$ induces a
perfect pairing of $O_K$-modules, up to a factor in~$\QQ^*$; and
2:~$X'$ is the set of Hodge structures of type $(-1,0),(0,-1)$ on the
$K$-vector space $K^2$ such that, up to sign, $\trace\circ\psi$ is a
polarization. The moduli space over $\QQ$ of triplets
$(A/S,\alpha,\lambda)$ with $S$ a $\QQ$-scheme, $\alpha\colon O_K\to
\End_S(A)$ a morphism of rings, and $\lambda$ a principal
$O_K$-polarization, is then the canonical model $S'_\QQ$ of~$S'_\CC$
(see also \cite[1.27]{Rapoport1} and \cite[4.11]{Deligne2}).

For $n\geq 1$, let $H_n$ be the kernel of the morphism $G'(\Zhat)\to
G'(\ZZ/n\ZZ)$, and let $S'_{\QQ,n}$ denote the Shimura variety
$\Sh_{H_n}(G',X')_\QQ$. Then $S'_{\QQ,n}$ is the moduli space for
4tuples $(A/S,\alpha,\lambda,\phi)$, with $S$ a $\QQ$-scheme,
$(A/S,\alpha,\lambda)$ an abelian scheme over $S$ with multiplications
by $O_K$ and a principal $O_K$-polarization, and with $\phi$ an
isomorphism of $S$-group schemes with $O_K$-action:
$$
\phi\colon (O_K/nO_K)_S^2\lto A[n], 
$$
such that there exists a (necessarily unique) isomorphism
$\ol{\phi}\colon (\ZZ/n\ZZ)_S \to \mu_{n,S}$ making the diagram:
$$
\begin{array}{ccc}
((O_K/nO_K)^2_S)^2 & \overset{\phi}{\lto} & A[n]^2\\
\downarrow\psi_n & & \downarrow e_{\lambda,n}\\
O_K\otimes(\ZZ/n\ZZ)_S & \overset{\id\otimes\ol{\phi}}\lto{} & 
O_K\otimes\mu_{n,S}
\end{array}
$$
commutative. In this diagram, $\psi_n$ is the pairing given by
$(\begin{smallmatrix}0 & 1\\-1 & 0\end{smallmatrix})$, and
$e_{\lambda,n}$ is the perfect pairing on $A[n]$ induced
by~$\lambda$. For $n\geq3$ the objects $((A/S,\alpha,\lambda))$ have
no non-trivial automorphisms (see \cite[IV, Thm.~5]{Mumford1}), and
$S'_{\QQ,n}$ is a fine moduli space.  (Representability by an
algebraic space can be found
in~\cite[\S1.23]{Rapoport1}. Quasi-projectiveness follows
from~\cite{BailyBorel1}.) In particular, for $n\geq3$ we do have a
polarized variation of $\ZZ$-Hodge structure on $S'_n(\CC)$, given by
the first homology groups of the fibers of the universal family.

\section{Monodromy and generic Mumford-Tate groups.}\label{sec3.4}
We recall some results that can be found in
\cite[Sections~1.1--1.3]{Moonen3}, with references to
\cite{CDK} and~\cite{Andre4}.

The Mumford-Tate group $\MT(V)$ of a $\QQ$-Hodge structure $V$, given
by $h\colon\SS\to\GL(V)_\RR$, is defined to be the smallest algebraic
subgroup $H$ of $\GL(V)_\QQ$ such that $h$ factors through~$H_\RR$.
Equivalently, $\MT(V)$ is the intersection in $\GL(V)$ of all
stabilizers of all lines generated by Hodge classes (i.e., of some
type $(p,p)$) in all $\QQ$-Hodge structures of the form
$\oplus_iV^{\otimes n_i}\otimes(V^*)^{\otimes m_i}$.

For $S$ a smooth complex algebraic variety with a polarizable
variation of $\QQ$-Hodge structure $V$ on the associated analytic
variety $S(\CC)$, there is a countable union $\Sigma$ of proper
algebraic subvarieties such that $s\mapsto \MT(V_s)$ is locally
constant outside $\Sigma$ (this makes sense because $V$ is a locally
constant sheaf on~$S(\CC)$). The smallest such $\Sigma$ is called the
Hodge exceptional locus, and its complement the Hodge generic
locus. For $s$ in $S(\CC)$ and not in $\Sigma$,
$\MT(V_s)\subset\GL(V_s)$ is called the generic Mumford-Tate group
(at~$s$).

Assume now that $S$ is connected, and that we have an element $s$
of~$S(\CC)$. Then the locally constant sheaf $V$ corresponds to a
representation $\rho\colon\pi_1(S(\CC),s)\to\GL(V_s)$, called the
monodromy representation. The algebraic monodromy group is defined to
be the smallest algebraic subgroup $H$ of $\GL(V_s)$ such that $\rho$
factors through $H$, i.e., it is the Zariski closure of the image
of~$\rho$; its connected component of identity is called the connected
algebraic monodromy group, and denoted $\AM(V_s)$. With these
hypotheses, we have the following theorem.
\begin{subthm}[Andr\'e]\label{thm3.4.1}
Assume moreover that $V$ admits a $\ZZ$-structure, that $s$ in
$S(\CC)$ is Hodge generic, and that there is a point $t$ in $S(\CC)$
such that $\MT(V_t)$ is abelian (i.e., $t$ is special). Then
$\AM(V_s)$ is the derived subgroup $\MT(V_s)^\der$ of $\MT(V_s)$,
i.e., the algebraic subgroup generated by commutators.
\end{subthm}
Let us now consider what this theorem implies for the variation of
Hodge structure that we have on $S'_n(\CC)$ ($n\geq3$), and, more
importantly, for its restrictions to subvarieties of~$S'_n(\CC)$.  The
Hodge exceptional locus of $S'_n(\CC)$ is by construction the union of
all lower dimensional subvarieties of Hodge type. The generic
Mumford-Tate group on $S'_n(\CC)$ is $G'$ (use that it contains a
subgroup of finite index of $G'(\ZZ)$, and that for all
$h=(h_1,h_2)\colon \CC^*\to\GL_2(\RR)^2$ in $X'$ one has
$\det(h_1(z))=z\bar{z}=\det(h_2(z))$ for all~$z$).
\begin{subprop}\label{prop3.4.2}
Let $n\geq3$. Let $C$ be an irreducible curve in $S'_{\CC,n}$ (i.e.,
an irreducible closed subvariety of dimension one); let $C^\nor$
denote its normalization and $C^\sm$ its smooth locus. Then $C$ is of
Hodge type if and only if the generic Mumford-Tate group on $C^\sm$ is
strictly smaller than~$G'_\QQ$. If $C$ is not of Hodge type and
contains a special point, then the connected algebraic monodromy group
on $C^\nor$ equals ${G'}^\der_\QQ = \Res_{K/\QQ}\SL_{2,K}$.
\end{subprop}
\begin{proof}
Suppose that $C$ is of Hodge type. Then some element in some tensor
construction of the variation of Hodge structure on $S'_{\CC,n}$ is a
Hodge class on $C$, but not on~$S'_{\CC,n}$. The interpretation of the
Mumford-Tate group as stabilizer of lines generated by Hodge classes
shows that the generic Mumford-Tate group on $C^\sm$ is strictly
smaller than~$G'_\QQ$.  Now suppose that the generic Mumford-Tate
group on $C^\sm$ is strictly smaller than~$G'_\QQ$. Then $C$ does
carry an extra Hodge class. The locus where this class is a Hodge
class is necessarily of dimension one, hence, $C$, being an
irreducible component of it, is of Hodge type. The second statement
follows now from Andr\'e's theorem above.
\end{proof}

\section{Irreducibility of images under Hecke correspondences.}
\label{sec3.5}
For $(G,X)$ a Shimura datum, $K_1$ and $K_2$ open subgroups of
$G(\AAf)$, and $g$ in $G(\AAf)$, one has the so-called Hecke
correspondence $T_g$ that is defined as follows. Consider the diagram:
$$
\Sh_{K_1}(G,X)_\CC \overset{\pi_1}{\longleftarrow} \Sh(G,X)_\CC 
\overset{{\cdot}g}{\longrightarrow}\Sh(G,X)_\CC 
\overset{\pi_2}{\longrightarrow}\Sh_{K_2}(G,X)_\CC, 
$$
where $\pi_1$ and $\pi_2$ are the quotient maps for the actions by
$K_1$ and $K_2$, respectively. The morphism $\pi_2\circ{\cdot}g$ is
the quotient for the action of $gK_2g^{-1}$, hence $\pi_1$ and
$\pi_2\circ{\cdot}g$ both factor through the quotient by $K:=K_1\cap
gK_2g^{-1}$, and $T_g$ is the correspondence:
$$
\Sh_{K_1}(G,X)_\CC \overset{\ol{\pi_1}}{\longleftarrow} \Sh_K(G,X)_\CC 
\overset{\ol{\pi_2\circ{\cdot}g}}{\longrightarrow}\Sh_{K_2}(G,X)_\CC. 
$$
Of course, $T_g$ exists already over the reflex field $E$ of $(G,X)$.
In particular, for $Z$ a closed subvariety of $\Sh_{K_1}(G,X)_E$, its
image $T_gZ$ is a closed subvariety of $\Sh_{K_2}(G,X)_E$.

We now specialize to our situation, i.e., to the Shimura datum
$(G',X')$ as above. For $p$ a prime number, we let $T_p$ be the Hecke
correspondence on $S'_\QQ$ given by the element $g(p)$ in $G'(\AAf)$
with $g(p)_p=(\begin{smallmatrix}p & 0\\0 & 1\end{smallmatrix})$ and
$g(p)_l=1$ for $l$ different from~$p$. Note that $g(p)^{-1}$ gives the
same correspondence as $g(p)$ does, because
$g(p)=(\begin{smallmatrix}p&0\\0&p\end{smallmatrix})
(\begin{smallmatrix}0&1\\1&0\end{smallmatrix})g(p)^{-1}
(\begin{smallmatrix}0&1\\1&0\end{smallmatrix})$.  The modular
interpretation of $T_p$ is the following. Let $[(A,\lambda)]$ in
$S'(\CC)$ denote the isomorphism class of a complex abelian surface
$A$ with multiplications by $O_K$ and with a principal
$O_K$-polarization~$\lambda$. Then, as a cycle, the image of
$[(A,\lambda)]$ is given by:
$$
T_p[(A,\lambda)] = \sum_H [(A/H,\ol{p\lambda})], 
$$
where $H$ ranges through the $O_K/pO_K$-submodules of $A[p](\CC)$ that
are free of rank one, and where $\ol{p\lambda}$ is the principal
$O_K$-polarization induced by $p\lambda$ on~$A/H$. In order to see
this, one uses, as in Section~\ref{sec3.3}, that $G'(\AAf)/G'(\Zhat)$
is the set of $O_K$-lattices in $K^2$ on which
$\psi=(\begin{smallmatrix}0&1\\-1&0\end{smallmatrix})$ induces a
perfect pairing of $O_K$-modules, up to a factor in $\QQ^*$, and that
the correspondence on it induced by $g(p)^{-1}$ sends such a lattice
to the set of lattices containing it with quotient free of rank one
over~$O_K/pO_K$.

\begin{subprop}\label{prop3.5.1}
Let $C$ be an irreducible curve in~$S'_\CC$. Suppose that $C$ is not
of Hodge type and that it contains a special point. Then, for all
primes $p$ large enough, $T_pC$ is irreducible.
\end{subprop}
\begin{proof}
Let $n\geq 3$ be some integer, and let $C_n$ be an irreducible
component of the inverse image of $C$ in~$S'_{\CC,n}$. Irreducibility
of $T_pC_n$ implies that of~$T_pC$. Let $V$ denote the polarized
variation of $\ZZ$-Hodge structure on $S'_n(\CC)$ that we considered
before, let $s$ be in~$C_n(\CC)$. We choose an isomorphism of
$O_K$-modules from $O_K^2$ to~$V_s$.  Let
$\rho\colon\pi_1(C_n(\CC),s)\to\SL_2(O_K)$ be the monodromy
representation. Proposition~\ref{prop3.4.2} implies that the Zariski
closure in $G'$ of $\rho(\pi_1(C_n^\nor(\CC),s))$ is the subgroup
$\Res_{O_K/\ZZ}\SL_{2,O_K}$. For $p$ prime, the correspondence $T_p$
on $S'_{\CC,n}$ is given by a diagram:
$$
S'_{\CC,n} \overset{\pi_1}{\longleftarrow} S'_{\CC,n,p} 
\overset{\pi_2}{\lto} S'_{\CC,n}. 
$$
For $T_pC_n=\pi_2\pi_1^{-1}C_n$ to be irreducible, it suffices that
$C_{n,p}$ be irreducible, with $C_{n,p}$ the covering of $C_n^\nor$
obtained from~$\pi_1$. But this covering corresponds to the
$\pi_1(C_n^\nor(\CC),s)$-set $\PP^1(O_K/pO_K)$ of
$O_K/pO_K$-submodules of $(O_K/pO_K)^2$ that are free of rank
one. Nori's Theorem \cite[Thm.~5.1]{Nori1} (Theorem~\ref{thm3.5.2}
below) implies that for $p$ large enough, the reduction map from
$\pi_1(C_n^\nor(\CC),s)$ to $\SL_2(O_K/pO_K)$ is surjective. Since
$\SL_2(O_K/pO_K)$ acts transitively on $\PP^1(O_K/pO_K)$,
irreducibility follows.
\end{proof}
\begin{subthm}[Nori]\label{thm3.5.2}
Let $\pi$ be a finitely generated subgroup of $\GL_n(\ZZ)$, let $H$ be
the Zariski closure of $\pi$, and for $p$ prime, let $\pi(p)$ be the
image of $\pi$ in~$\GL_n(\FF_p)$. Then, for almost all $p$, $\pi(p)$
contains the subgroup of $H(\FF_p)$ that is generated by the elements
of order~$p$.
\end{subthm}

\section{Galois action.}\label{sec3.6}
The aim of this section is to show that the Galois orbits of special
points in $S'(\Qbar)$ are big, in a suitable sense. For $A$ and $B$
abelian surfaces (over some field) with $O_K$-action, we let
$\Hom_{O_K}(A,B)$ be the $O_K$-module of morphisms from $A$ to $B$
that are compatible with the $O_K$-actions.
\begin{sublem}
Let $x$ in $S'(\Qbar)$ be a special point, corresponding to a triplet
$(A,\alpha,\lambda)$ with $A$ an abelian surface over $\Qbar$,
$\alpha\colon O_K\to\End(A)$ and $\lambda$ a principal
$O_K$-polarization.  Then $\End_{O_K}(A)$ is an order, containing
$O_K$, of a totally imaginary quadratic extension of~$K$.
\end{sublem}
\begin{proof}
Let $R$ be the endomorphism algebra $\QQ\otimes\End(A)$ of~$A$. Then
$R$ is a semi-simple $\QQ$-algebra containing a commutative
semi-simple subalgebra of dimension~$4$. Suppose that $A$ is
simple. Then $R$ is a division algebra.  Since $R$ acts faithfully on
$\rH_1(A(\CC),\QQ)$, it has dimension dividing $4$, hence $R$ is a
quadratic extension of~$K$. Since $\RR\otimes R$ has a complex
structure commuting with the $R$-action, $R$ is a totally
imaginary. Suppose now that $A$ is not simple. Then $A$ is isogeneous
to the product of two elliptic curves, $B_1$ and $B_2$, say. These
elliptic curves are in fact isogeneous to each other, because
otherwise $K$ does not admit a morphism to the endomorphism algebra of
$B_1\times B_2$. So $A$ is isogeneous to $B^2$, with $B$ some elliptic
curve. Since $A$ is of CM-type, $\QQ\otimes\End(B)$ is an imaginary
quadratic field $E$, and $R=\rM_2(E)$.  In this case $\End_{O_K}(A)$
is an order in the totally imaginary extension $K\otimes E$ of~$K$.
\end{proof}

\begin{subthm}\label{thm3.6.2}
There exist real numbers $\eps>0$ and $c>0$ such that for
$(A,\alpha,\lambda)$ corresponding to a special point $x$ in
$S'(\Qbar)$ one has:
$$
|\Gal(\Qbar/\QQ){\cdot}x| >c\,|\discr(R_x)|^\eps, 
$$
where $R_x=\End_{O_K}(A)$. 
\end{subthm}
\begin{rem}
The proof will show that one can take $\eps$ to be any number less
than~$1/4$. (To get this, one also has to optimize
Theorem~\ref{thm3.6.4}, noting that we only apply Stark's result to
fields $L$ of degree at least~4.)  Assuming the generalized Riemann
hypothesis at this point does not improve this exponent (this is
caused by the case where $\QQ\otimes R_x$ is Galois over $\QQ$ with
group $(\ZZ/2\ZZ)^2$).
\end{rem}
\begin{proof}
Let $f\colon S'_\QQ\to S_\QQ$ be the morphism induced by the closed
immersion of the Shimura data $(G_\QQ',X')\to(G_\QQ,X)$. Since $f$ is
finite, and since the Hecke correspondences on $S_\QQ$ permute the
irreducible components transitively, the statement we want to prove is
equivalent to its analog for~$S_\QQ$. So we will show in fact that
there are positive $\eps$ and $c$ such that for $x$ special in
$S(\Qbar)$ corresponding to $(A,\alpha)$, we have:
$$
|\Gal(\Qbar/\QQ){\cdot}x| >c\,|\discr(R_x)|^\eps. 
$$
For $x$ special in $S(\Qbar)$, let $L_x$ be $\QQ\otimes R_x$, and let
$M_x$ be the Galois closure in $\Qbar$ of~$L_x$. Since $M_x$ is of
degree at most $8$ over $\QQ$, the statement we want to prove is
equivalent to the existence of positive $\eps$ and $c$ such that for
all special $x$ in~$S(\Qbar)$:
$$
|\Gal(\Qbar/M_x){\cdot}x| >c\,|\discr(R_x)|^\eps. 
$$
So let now $x$ be special in $S(\Qbar)$, corresponding to some
$(A_x,\alpha_x)$.  To study the $\Gal(\Qbar/M_x)$-orbit of $x$, we
construct a zero-dimensional subvariety of Hodge type over $M_x$,
containing $x$, and we use the theory of Shimura-Taniyama on complex
multiplication, rephrased in the language of Shimura varieties (see
\cite[Section~2.2]{Moonen2}).  Let $H_x$ be $\rH_1(A(\CC),\ZZ)$; it is
an $R_x$-module that is locally free of rank two as $O_K$-module. Let
$f\colon K^2\to \QQ\otimes H_x$ be an isomorphism of $K$-vector
spaces.  The Hodge structure on $H_x$ given by $A_x$ gives an element
$h_x$ of~$X$.  The lattice $f^{-1}H_x$ in $K^2$ corresponds to an
element $\ol{g_x}$ of $G(\AAf)/G(\Zhat)$. By construction, $x$ is the
image of~$(h_x,g_x)$.  Let $T_x:=\Res_{L_x/\QQ}\Gm[,L_x]$. Then $f$
gives a closed immersion $T_x\to G_\QQ$. Since $\QQ\otimes H_x$ is a
one-dimensional $L_x$-vector space, $T_x$ is its own centralizer
in~$G_\QQ$. It follows that $h_x$ factors through~$T_{x,\RR}$. Hence
we have a closed immersion of Shimura data:
$(T_x,\{h_x\})\to(G_\QQ,X)$. The reflex field of $(T_x,\{h_x\})$ is
contained in $M_x$, hence we have a canonical model
$\Sh(T_x,\{h_x\})_{M_x}$ over~$M_x$.  We put $U_x:=T_x(\AAf)\cap
g_xG(\Zhat)g_x^{-1}$. Then one easily verifies that we have an
injective morphism of Shimura varieties
$\Sh_{U_x}(T_x,\{h_x\})_{M_x}\to S_{M_x}$, which, on $\CC$-valued
points, is given by $\ol{t}\mapsto\ol{(h_x,tg_x)}$. By construction,
$U_x$ is the stabilizer in $T_x(\AAf)$ of the lattice $f^{-1}H_x$; it
follows that $U_x=R_x^{\swedge,*}$, hence:
$$
{\Sh_{U_x}(T_x,\{h_x\})_{M_x}}(\Qbar)=
L_x^*\backslash(\AAf\otimes L_x)^*/R_x^{\swedge,*} = \Pic(R_x). 
$$ 
Our next objective is to describe in sufficient detail the action of
$\Gal(\Qbar/M_x)$ on $\Pic(R_x)$ induced by the above bijections.
Class field theory gives a continuous surjection from
$M_x^*\backslash\AAf[M_x,]^*$ to $\Gal(\Qbar/M_x)^\ab$, characterized
by the following  property. In a representation of $\Gal(\Qbar/M_x)^\ab$
that is unramified at a finite place $v$ of $M_x$, the arithmetic
Frobenius element is the image of the class of an id\`ele that is
trivial at all places other than $v$, and the inverse of a uniformizer
at~$v$.  Let $\mu\colon\Gm[,\CC]\to\SS_\CC$ be the cocharacter obtained
by composing $\CC^*\to \CC^*\times\CC^*$, $z\mapsto(z,1)$ with the
inverse of the isomorphism
$\SS(\CC)=(\CC\otimes_\RR\CC)^*\to\CC^*\times\CC^*$, $x\otimes
y\mapsto(xy,x\ol{y})$. Then $h_x\circ\mu$ is defined over $M_x$, and
one defines:
$$
r_x\colon T'_x:=\Res_{M_x/\QQ}\Gm[,M_x]\to T_x
$$
to be the morphism $\Res_{M_x/\QQ}(h_x\circ\mu)$ composed with the
norm map from $\Res_{M_x/\QQ}T_{x,M_x}$ to~$T_x$. With these
definitions, the quotient $\Gal(\Qbar/M_x)^\ab$ of $T'_x(\AAf)$ acts
on $\Pic(R_x)$ via the morphism $r_x$, where we view $\Pic(R_x)$ as
$T_x(\QQ)\backslash T_x(\AAf)/R_x^{\swedge,*}$. It follows that:
$$
|\Gal(\Qbar/M_x){\cdot}x| = 
|\text{image of $r_x(T_x'(\AAf))$ in $\Pic(R_x)$}|. 
$$

We will need a more explicit description of $r_x$, in terms of the CM
type associated to~$h_x$. The morphism $h_x\colon \CC^*\to(\RR\otimes
L_x)^*$ extends to a morphism of $\RR$-algebras
$h\colon\CC\to\RR\otimes L_x$.  Extending scalars from $\RR$ to $\CC$
gives a morphism of $\CC$-algebras $\id\otimes
h\colon\CC\otimes_\RR\CC\to\CC\otimes L_x$. Via the isomorphisms:
$$
\CC\otimes_\RR\CC\to\CC\times\CC,\quad x\otimes y\mapsto(xy,x\ol{y}),
$$ 
and 
$$
\CC\otimes L_x\to\CC^{\Hom(L_x,\CC)},\quad x\otimes y\mapsto 
(\phi\mapsto x\phi(y)), 
$$
the idempotent $(1,0)$ of $\CC\times\CC$ gives an idempotent in
$\CC^{\Hom(L_x,\CC)}$, i.e., a partition of $\Hom(L_x,\CC)$ into two
sets $\Phi_x$ and $\iota\Phi_x$, where $\iota$ is the complex
conjugation on~$\CC$. The set $\Phi_x$ is the CM type corresponding
to~$h_x$. Since $M_x$ is the Galois closure of $L_x$ in $\CC$,
$\Hom(L_x,M_x)=\Hom(L_x,\CC)$. With these notations, we have, for any
$\QQ$-algebra~$R$:
$$
r_x\colon (R\otimes M_x)^*\lto \prod_{\phi\in\Phi_x}(R\otimes
L_x)^*,\quad u\longmapsto \prod_{\phi\in\Phi_x}\Norm_\phi(u),
$$
where $\Norm_\phi$ is the norm map of the extension $\phi\colon
R\otimes L_x\to R\otimes M_x$.  Finally, let $\phi_0$ be in
$\Hom(L_x,M_x)$, and define
$\Sigma_{x,\phi_0}:=\{g\in\Gal(M_x/\QQ)\,|\,g\phi_0\in\Phi_x\}$. Then
we have:
$$
\phi_0\circ r_x\colon T'_x\lto T_x\into T'_x, \quad 
u\mapsto \prod_{g\in\Sigma_{x,\phi_0}}g^{-1}u, 
$$
for all $\QQ$-algebras $R$ and all $u$ in $(R\otimes M_x)^*$. This is
the description of $r_x$ that we work with.

Since $M_x$ is generated over $K$ by the extension $L_x$ and its
conjugate, $M_x$ has degree 4 or 8 over $\QQ$, and its Galois group
$\Gal(M_x/\QQ)$ is isomorphic to $\ZZ/4\ZZ$, $\ZZ/2\ZZ\times\ZZ/2\ZZ$,
or $D_4$, the dihedral group of order~8. We define $T$ to be
$\Res_{K/\QQ}\Gm[,K]$; note that $T$ is a subtorus of $T_x$, equal to
the center of $G_\QQ$. We will see below that $r_x\circ\phi_0\colon
T_x\to T_x$ induces an endomorphism of $T_x/T$ whose image, after
passing to $\AAf$-valued points, in $\Pic(R_x)/\Pic(O_K)$ is big
enough for our purposes.

Suppose first that $\Gal(M_x/\QQ)$ is isomorphic to $\ZZ/4\ZZ$, say
with generator~$\sigma$. Then $M_x=L_x$, $\sigma^2$ is the complex
conjugation and $K=L_x^{\ld\sigma^2\rd}$. After changing $\phi_0$, if
necessary, one has that $\Sigma_{x,\phi}=\{1,\sigma\}$. The formula
above for $\phi_0\circ r_x$ shows that $r_x$ is simply given by the
element $1+\sigma^{-1}$ of $\ZZ[\Gal(M_x/\QQ)]$. Since $\sigma^2$ acts
as $-1$ on $T_x/T$, we have $r_x\circ(1-\sigma^{-1})=2$ on~$T_x/T$. It
follows that:
$$
|\Gal(\Qbar/M_x){\cdot}x| \geq 
|\text{image of ${\cdot}2\colon \Pic(R_x)/\Pic(O_K)\to\Pic(R_x)/\Pic(O_K)$}|. 
$$
Theorem~\ref{thm3.6.4} below finishes the proof in this case. 

Suppose now that $\Gal(M_x/\QQ)$ is isomorphic to
$\ZZ/2\ZZ\times\ZZ/2\ZZ$.  After changing $\phi_0$, if necessary, one
has $\Sigma_{x,\phi_0}=\{1,\sigma\}$, with $\sigma$ of order two and
$K_x:=L_x^{\ld\sigma\rd}\neq K$.  Let $R_x'$ be the order $O_{K_x}\cap
R_x$ of~$K_x$. Since $r_x$ is given by $1+\sigma$, the induced map
$T_x'(\AAf)\to\Pic(R_x)$ factors through $\Pic(R_x')\to\Pic(R_x)$
induced by the inclusion $R_x'\to R_x$.  The fact that $\sigma$ acts
as $1$ on $T_x'':=\Res_{K_x/\QQ}\Gm[,K_x]$ and as $-1$ on $T_x/T_x''$
implies that the kernel of the map $\Pic(R_x')\to\Pic(R_x)$ is killed
by multiplication by~$2$.  Since $1+\sigma$ acts as multiplication by
$2$ on $T_x''$, we get:
$$
|\Gal(\Qbar/M_x){\cdot}x| \geq 
|\text{image of ${\cdot}4\colon \Pic(R_x')\to\Pic(R_x')$}|. 
$$
Since the order $O_K\otimes R_x'$ of $L_x$ is contained in $R_x$, 
and has discriminant $\discr(R_x')^2\discr(O_K)^2$, we have: 
$$
|\discr(R_x')|\geq |\discr(O_K)|^{-1}|\discr(R_x)|^{1/2}. 
$$
The proof in this case is finished by Theorem~\ref{thm3.6.4}. 

Suppose that $\Gal(M_x/\QQ)$ is isomorphic to~$D_4$. Let $\tau$ and
$\sigma$ be generators of $\Gal(M_x/\QQ)$, with
$M_x^{\ld\tau\rd}=L_x$, and with $\sigma$ of order~$4$. Then
$\sigma^2$ is the complex conjugation, and
$\tau\sigma^{-1}=\sigma\tau$. After changing $\phi_0$, if necessary,
we have $\Sigma_{x,\phi_0}=\{1,\tau,\sigma,\sigma\tau\}$. It follows
that $\phi_0\circ r_x$ is given by the element
$t:=1+\tau+\sigma^3+\sigma\tau$ of $\ZZ[\Gal(M_x/\QQ)]$. Using that
$\Norm_{\phi_0}=1+\tau$, a simple computation gives:
$$
\phi_0\circ r_x\circ\phi_0\circ\Norm_{\phi_0} = 
2(1+\tau) + \sigma\Norm_{M_x/K}. 
$$
It follows that $r_x\circ\phi_0$ acts as $2$ on $T_x/T$. We conclude that: 
$$
|\Gal(\Qbar/M_x){\cdot}x| \geq 
|\text{image of ${\cdot}2\colon \Pic(R_x)/\Pic(O_K)\to\Pic(R_x)/\Pic(O_K)$}|. 
$$
Theorem~\ref{thm3.6.4} finishes the proof in this last case. 
\end{proof}

\begin{subthm}\label{thm3.6.4}
Let $K$ be a totally real number field. There exists $c>0$ such that
for all orders $R$, containing $O_K$, in totally complex quadratic
extensions $L$ of $K$, one has:
$$
|\text{image of ${\cdot}4$ on $\Pic(R)/\Pic(O_K)$}| \geq c\,|\discr(R)|^{1/8}. 
$$
If one assumes the generalized Riemann hypothesis, then one can
replace the exponent $1/8$ by any number less than~$1/2$.
\end{subthm}
\begin{proof}
We will use the following lower bound for class numbers: 
\begin{quote}\sl 
let $K$ be a totally real number field; there exists $c>0$ such that
for all totally complex quadratic extensions $L$ of $K$, one has:
$$
|\Pic(O_L)| \geq c\,|\discr(O_L)|^{1/6}. 
$$
\end{quote}
In order to prove this, one distinguishes two cases: $K=\QQ$ and
$K\neq\QQ$, and one notes that the regulator $\Reg(O_L)$ is at
most~$\Reg(O_K)$. In the case $K\neq\QQ$ one uses the following
consequence of Stark's Theorem~2 in \cite{Stark1}: 
\begin{quote}\sl 
for $K$ a totally real number field, there exists $c>0$ such that for
all totally complex quadratic extensions $L$ of $K$, one has:
$$
|\Pic(O_L)|\geq c\,|\discr(O_L)|^{1/2-1/[L:\QQ]}. 
$$ 
\end{quote}
In the case $K=\QQ$ one applies the Brauer-Siegel theorem (see for
example \cite[Ch.~XVI]{Lang1}):
\begin{quote} \sl
for $N>0$ and $\eps>0$, there exists $c>0$ such that for all Galois 
extensions $L$ of $\QQ$ of degree at most $N$ one has: 
$$
|\Pic(O_L)|{\cdot}\Reg(O_L) \geq c\,|\discr(O_L)|^{1/2-\eps}.
$$
\end{quote}
Combining these two results, and using that $[L:\QQ]\geq3$ if
$K\neq\QQ$ gives the inequality we want. We could replace the exponent
$1/6$ by $1/4$ if we would just use that $[L:\QQ]\geq4$ if $K\neq\QQ$.

Let now $K$, $R$ and $L$ be as in the theorem. Then $R$ is the inverse
image of a subring $\ol{R}$ of some finite quotient $\ol{O_L}$
of~$O_L$.  We have an exact sequence:
$$
0\lto R^*\lto O_L^*\lto \ol{O_L}^*/\ol{R}^* \lto\Pic(R)\lto\Pic(O_L)\lto 0. 
$$
The torsion of $O_L^*$ is bounded in terms of the degree of $K$, and
by Dirichlet's theorem on units the quotient $O_L^*/O_K^*$ is finite.
The long exact cohomology sequence obtained by taking
$\Gal(L/K)$-invariants of the short exact sequence: 
$$
0\lto \tors(O_L^*)\lto O_L^* \lto O_L^*/\tors(O_L^*)\lto 0
$$
gives an injection from $(O_L^*/\tors(O_L^*))/(O_K^*/\tors(O_K^*))$
into $\rH^1(\Gal(L/K),\tors(O_L^*))$, showing that
$(O_L^*/\tors(O_L^*))/(O_K^*/\tors(O_K^*))$ is of order at most
two. We conclude that there exists $c>0$, depending only on the degree
of $K$, such that:
$$
|\Pic(R)| \geq c\,\frac{|\ol{O_L}^*|}{|\ol{R}^*|}\,|\Pic(O_L)|. 
$$
On the other hand, we have: 
$$
\discr(R) = \left(\frac{|\ol{O_L}|}{|\ol{R}|}\right)^2\discr(O_L). 
$$
We claim that for every $\eps>0$ there exists $c>0$, depending only on 
the degree of $K$, such that: 
$$
\frac{|\ol{O_L}^*|}{|\ol{R}^*|} \geq 
c\left(\frac{|\ol{O_L}|}{|\ol{R}|}\right)^{1-\eps}. 
$$
To prove this claim, one notes that: 
$$
\frac{|\ol{R}|}{|\ol{R}^*|} = \prod_{\text{$k$ res field of
$\ol{R}$}}\frac{|k|}{|k^*|}, \quad \text{and}\quad 
\frac{|\ol{O_L}|}{|\ol{O_L}^*|} = \prod_{\text{$k$ res field of
$\ol{O_L}$}}\frac{|k|}{|k^*|}. 
$$
A simple computation then shows:
$$
\frac{|\ol{O_L}^*|}{|\ol{R}^*|}\geq 
n\prod_{p|n}\left(1-\frac{1}{p}\right)^{[L:\QQ]}
\geq n\left(\frac{1}{5\log(n)}\right)^{[L:\QQ]}, 
$$
where $n=|\ol{O_L}|/|\ol{R}|$ is assumed to be at least~2. 
We conclude that there exists $c>0$, depending only on $K$, such that: 
$$
|\Pic(R)| \geq c\,|\discr(R)|^{1/7}. 
$$
In order to finish the proof of the theorem, it suffices to prove that
for every $\eps>0$ there exists $c>0$, depending only on $K$, such
that $|\Pic(R)[2]|\leq c\,|\discr(R)|^\eps$. To do this, we proceed in
the same way as we did in~\cite[Lemma~3.4]{Edixhoven1}. As $\Pic(R)$
is a finite commutative group, the two $\FF_2$-vector spaces
$\Pic(R)[2]$ and $\FF_2\otimes\Pic(R)$ have the same dimension.  The
cover of $\Spec(R)$ by the disjoint union of $\Spec(\ZZ_2\otimes R)$
and $\Spec(R[1/2])$ gives an exact sequence:
$$
(\QQ_2\otimes L)^* \lto \Pic(R) \lto \Pic(R[1/2]) \lto 0. 
$$
It follows that $\dim_{\FF_2}\FF_2\otimes\Pic(R)$ is bounded by
$\dim_{\FF_2}\FF_2\otimes\Pic(R[1/2])$ plus a number depending only on
the degree of~$K$. We put $S:=\Spec(R[1/2])$ and
$T:=\Spec(O_K[1/2])$. The long exact sequence associated to the
multiplication by two on the sheaf $\Gm$ on the etale site $S_\et$ of
$S$ shows that $\dim_{\FF_2}\FF_2\otimes\Pic(S)$ is at most
$\dim_{\FF_2}\rH^1(S_\et,\FF_2)$. Let $\pi\colon S\to T$ be the
morphism induced by the inclusion of $O_K$ in~$R$. Then
$\rH^1(S_\et,\FF_2)$ is the same as $\rH^1(T_\et,\pi_*\FF_2)$, and we
have a short exact sequence:
$$
0\lto \FF_{2,S_\et} \lto \pi_*\pi^*\FF_{2,T} \lto j_!\FF_{2,U} \lto 0, 
$$
where $j\colon U\to T$ is the maximal open immersion over which $\pi$
is etale. Let $i\colon Z\to T$ be closed immersion giving the
complement of $U$, with $Z$ reduced. Then the long exact sequences of
cohomology groups associated to the exact sequence above and to the
exact sequence:
$$
0\lto j_!\FF_{2,U_\et} \lto \FF_{2,S_\et} \lto i_*\FF_{2,Z_\et} \lto 0
$$
show that there exists an integer $c$, depending only on $K$, such that: 
$$
\dim_{\FF_2}\FF_2\otimes\Pic(R) \leq 
c+[K:\QQ]\,\left|\{\text{$p$ prime dividing $\discr(R)$}\}\right|. 
$$
As $2^{|\{p|n\}|}=n^{o(1)}$, for $n\to\infty$, we have proved the
first statement of the theorem. If one assumes GRH, then the
Brauer-Siegel theorem as stated above is true without the condition
that the extension $\QQ\to L$ be Galois, see \cite[XIII, \S4]{Lang1}.
\end{proof}

\section{Intersection numbers.}\label{sec3.7}
The aim of this section is to give a bound on intersections of
subvarieties of Shimura varieties, provided that they are finite. In
particular, we need to study the intersection of a subvariety with its
images under Hecke correspondences. As our arguments work for general
Shimura varieties, we give such a result in the general case. The main
tool used in proving the result is the Baily-Borel compactification,
together with its given ample line bundles. We start by recalling some
properties of these line bundles, that follow immediately from the
results in \cite{BailyBorel1} (see also \cite{BailyBorel2}).

\begin{subthm}\label{thm7.1}
Let $(G,X)$ be a Shimura datum. For $K\subset G(\AAf)$ a compact open
subgroup, let $S_K:=\Sh_K(G,X)_\CC$ the corresponding complex Shimura
variety, and $\ol{S}_K$ its Baily-Borel compactification. For every
such $K$, and for every sufficiently divisible positive integer $n$,
the $n$th power of the line bundle of holomorphic forms of maximal
degree of $X$ descends to $S_K$, and extends uniquely to a very ample
line bundle $\calL_{K,n}$ on $\ol{S}_K$, such that, at the generic
points of the boundary components of codimension one, it is given by
$n$th powers of forms with logarithmic poles. Let $K_1$ and $K_2$ be
compact open subgroups of $G(\AAf)$, and $g$ in $G(\AAf)$ such that
$K_2\subset gK_1g^{-1}$. Then the morphism from $S_{K_2}$ to $S_{K_1}$
induced by $g$ extends to a morphism $f\colon
\ol{S}_{K_2}\to\ol{S}_{K_1}$. If $n$ is positive and sufficiently
divisible so that $\calL_{K_1,n}$ exists, then $\calL_{K_2,n}$ exists,
and is canonically isomorphic to~$f^*\calL_{K_1,n}$.
\end{subthm}
\begin{proof}
Let us briefly recall how the compactification $\ol{S}_K$ is defined.
Let $X^+$ be a connected component of~$X$. Then each connected
component of $S_K(\CC)$ is of the form $\Gamma_i\backslash X^+$, with
$\Gamma_i$ an arithmetic subgroup of $G^\ad(\QQ)$ ($G^\ad$ being the
quotient of $G$ by its center). The compactification $\ol{S}_K(\CC)$
is then defined as the disjoint union of the
$\Gamma_i\backslash\ol{X^+}$, where $\ol{X^+}$ is the union of $X^+$
with its so-called rational boundary components, endowed with the
Satake topology. It follows that we can write $\ol{S}_K(\CC)$ as
$G(\QQ)\backslash(\ol{X}\times G(\AAf)/K)$, with $\ol{X}$ the disjoint
union of the~$\ol{X^+}$.

Let $X^\ad$ be the $G^\ad(\RR)$-conjugacy class of morphisms from
$\SS$ to $G^\ad_\RR$ containing the image of~$X$. Each connected
component of $X$ maps isomorphically to one of $X^\ad$ (see
\cite[1.6.7]{Moonen2}). We first prove the Theorem above for the
Shimura datum $(G^\ad,X^\ad)$. The group $G^\ad$ is a product of
simple algebraic groups $G_j$ over $\QQ$, and $X^\ad$ decomposes as a
product of~$X_j$'s. For compact open subgroups $K$, $K_1$ and $K_2$
that are products of compact open subgroups of the $G_j(\AAf)$, the
corresponding Shimura varieties decompose as a product, so that it
suffices to treat the $G_j$ separately. If $(G_j,X_j)$ gives compact
Shimura varieties, Kodaira's theorem
(\cite[Section~1.4]{GriffithsHarris}) implies what we want, for
compact open subgroups $K_j$ that are sufficiently small; for
arbitrary $K_j$ one takes quotients by finite groups. Suppose now that
$(G_j,X_j)$ does give Shimura varieties that are not compact. If $G_j$
is of dimension 3, then it is isomorphic to $\PGL_{2,\QQ}$, and we are
in the case of modular curves, where the Theorem we are proving is
well known (the canonical line bundle with log poles at the cusps on
the modular curve $X(n)$, $n\geq3$, has degree $>0$). Suppose now that
$G_j$ has dimension~$>3$. Then the boundary components are of
codimension $>1$, and the results we want are given in
\cite[Thm.~10.11]{BailyBorel1}.

The case of arbitrary open compact subgroups of $G^\ad(\AAf)$ follows
by considering quotients by finite groups. The theorem for $(G,X)$
itself follows from the fact that the connected components of the
$S_K(\CC)$ are of the $\Gamma\backslash X^+$, with $\Gamma$ an
arithmetic subgroup of~$G^\ad(\QQ)$.
\end{proof}

\begin{subthm}\label{thm7.2}
Let $(G,X)$ be a Shimura datum, let $K_1$ and $K_2$ be compact open
subgroups of $G(\AAf)$, and let $Z_1$ and $Z_2$ be closed subvarieties
of the Shimura varieties $S_1:=\Sh_{K_1}(G,X)_\CC$ and
$S_2:=\Sh_{K_2}(G,X)_\CC$, respectively. Suppose that $Z_1$ or $Z_2$
is of dimension at most one. Then there exists an integer
$c$ such that for all $g$ in $G(\AAf)$ for which $T_gZ_1\cap Z_2$ is
finite, one has: 
$$
|T_gZ_1\cap Z_2| \leq c\,\deg(\ol{\pi_1}\colon S_g\to S_1), 
$$
where $S_g=\Sh_{K_g}(G,X)_\CC$ with $K_g=K_1\cap gK_2g^{-1}$, and with
$T_g$ and $\ol{\pi_1}$ as in the beginning of Section~\ref{sec3.5}.
\end{subthm}
\begin{proof}
We start with two reductions. First of all, writing $Z_1$ and $Z_2$ as
the unions of their irreducible components, one sees that we may
suppose that $Z_1$ and $Z_2$ are irreducible. Secondly, for $g$ in
$G(\AAf)$, let $p_{1,g}$ and $p_{2,g}$ be the morphisms from $S_g$ to
$S_1$ and $S_2$, respectively. Then one has:
$$
T_gZ_1\cap Z_2 = p_{2,g}\left(p_{1,g}^{-1}Z_1\cap
p_{2,g}^{-1}Z_2\right), 
$$
which shows that $T_gZ_1\cap Z_2$ is finite if and only if
$p_{1,g}^{-1}Z_1\cap p_{2,g}^{-1}Z_2$ is, and that $|T_gZ_1\cap Z_2|$
is at most $|p_{1,g}^{-1}Z_1\cap p_{2,g}^{-1}Z_2|$. This also shows
that we may replace $K_1$ and $K_2$ by smaller compact open
subgroups. Hence we may suppose, by the previous theorem, that we have
very ample line bundles $\calL_1$ and $\calL_2$ on the Baily-Borel
compactifications $\ol{S_1}$ and $\ol{S_2}$ such that, for each $g$,
$\ol{p_{1,g}}^*\calL_1$ and $\ol{p_{2,g}}^*\calL_2$ are isomorphic to
the same line bundle $\calL_g$ on~$\ol{S_g}$. 

We let $\ol{Z_1}$ and $\ol{Z_2}$ be the closures of $Z_1$ and $Z_2$ in
$\ol{S_1}$ and $\ol{S_2}$, respectively. Let $m$ denote the degree of
$\ol{Z_2}$ with respect to~$\calL_2$. Let $g$ be in $G(\AAf)$, such
that $T_gZ_1\cap Z_2$ is finite. If the intersection is empty, there
is nothing to prove, so we suppose that the intersection is not
empty. Then the codimension of $Z_2$ is at least the dimension $d$ of
$Z_1$, and we can choose $f_1,\ldots,f_d$ in
$\rH^0(\ol{S_2},\calL_2^{\otimes m})$ such that $\ol{Z_2}$ is
contained in $V_{S_2}(f_1,\ldots,f_d)$, and $\ol{T_gZ_1}\cap
V_{S_2}(f_1,\ldots,f_d)$ is finite (because of our assumption on the
dimensions of $Z_1$ and $Z_2$, $\ol{Z_1}\cap\ol{Z_2}$ is finite). It
then follows that $|\ol{p_{1,g}^{-1}Z_1}\cap\ol{p_{2,g}^{-1}Z_2}|$ is
at most $m^d$ times the degree of $\ol{p_{1,g}^{-1}Z_1}$ with respect
to $\calL_g$. But this degree is $\deg(p_{1,g})$ times the degree of
$\ol{Z_1}$ with respect to $\calL_1$, hence we have: 
$$
|T_gZ_1\cap Z_2|\leq \deg(p_{1,g})m^d\deg_{\calL_1}(\ol{Z_1}). 
$$
\end{proof}

\section{Proof of the main results.}\label{sec3.8}
We will now prove Theorems~\ref{thm3.1} and~\ref{thm3.2}. We first deal with 
Thm.~\ref{thm3.1}. As we have already noticed, we may as well replace 
$S_\CC$ by $S'_\CC$, so let $C$ be an irreducible closed curve in $S'_\CC$ 
that contains infinitely many CM points. We have to show that $C$ is of  
Hodge type. 

Since $C$ has infinitely many points in $S'(\Qbar)$, it is, as a reduced 
closed subscheme of $S'_\CC$, defined over~$\Qbar$. To be precise, there 
is a unique closed subscheme $C_\Qbar$ of $S'_\Qbar$ that gives $C$ after 
base change from $\Qbar$ to~$\CC$. But then $C_\Qbar$ has only finitely many 
conjugates under $\Gal(\Qbar/\QQ)$; we let $C_\QQ$ be the reduced closed 
subscheme of $S'_\QQ$ that, after base change to $\Qbar$, gives the union 
of these conjugates. In other words, we simply let $C_\QQ$ be the image 
of $C$ under the morphism of schemes $S'_\CC\to S'_\QQ$. 

Let $x$ in $C(\Qbar)$ be a CM point, corresponding to 
a pair $(A,\lambda)$ with $A$ an abelian variety over $\Qbar$ with 
multiplications by $O_K$ and with $\lambda$ a principal $O_K$-polarization. 
As before, we let $R_x$ denote $\End_{O_K}(A)$, $L_x:=\QQ\otimes R_x$ and 
$M_x$ the Galois closure of $L_x$ in~$\CC$. Let $\ol{x}$ be the image of $x$ 
in~$S(\Qbar)$. In the proof of Theorem~\ref{thm3.6.2} we have seen that the 
quotient $\Gal(\Qbar/M_x)^\ab$ of $(\AAf\otimes M_x)^*$ acts on 
the subset $\Gal(\Qbar/M_x){\cdot}\ol{x}$ of 
$L_x^*\backslash(\AAf\otimes L_x)^*/R_x^{\swedge,*}$ via the morphism 
$$
r_x\colon (\AAf\otimes M_x)^*\lto (\AAf\otimes L_x)^*,\quad 
u\mapsto \prod_{\phi\in\Phi_x}\Norm_\phi(u),
$$ 
where the subset $\Phi_x$ of $\Hom(L_x,M_x)$ is the CM type
of~$x$. Since $\Phi_x$ is a set of representatives for the action of
$\Gal(L_x/K)$ on $\Hom(L_x,M_x)$, it follows that this map $r_x$
factors through the subgroup of elements of $(\AAf\otimes L_x)^*$
whose norm to $(\AAf\otimes K)^*$ is in~$\AAf^*$. Hence, in the
notation of the proof of Theorem~\ref{thm3.6.2}, the morphism $r_x$
factors through the intersection of the subtorus $T_x$ of $G_\QQ$
and~$G'_\QQ$. It follows that the action of $\Gal(\Qbar/M_x)$ on
$\Gal(\Qbar/M_x){\cdot}x$ is given by the same morphism $r_x$, taking
values in~$G'(\AAf)$. 
\begin{sublem} 
Suppose that $p$ is a prime that is split in $R_x$, i.e., such that
$\FF_p\otimes R_x$ is isomorphic to a product of copies
of\/~$\FF_p$. Then $\Gal(\Qbar/\QQ)x$ is contained in $C_\QQ(\Qbar)\cap
(T_pC_\QQ)(\Qbar)$.
\end{sublem}
\begin{proof}
The localization $\ZZ_{(p)}\otimes R_x$ of $R_x$ is the same as that
of~$O_K$.  Hence, if we let $H_x$ denote $\rH_1(A(\CC),\ZZ)$, then
$\ZZ_{(p)}\otimes H_x$ is free of rank one over $\ZZ_{(p)}\otimes
R_x$. It follows that we can choose the isomorphism $f\colon
K^2\to\QQ\otimes H_x$ to preserve the integral structures on both
sides at $p$, i.e., such that it induces an isomorphism from
$(\ZZ_{(p)}\otimes R_x)^2$ to $\ZZ_{(p)}\otimes H_x$.  We note that
$p$ is split in $M_x$ (i.e., $\FF_p\otimes O_{M_x}$ is a product of
copies of $\FF_p$), because $M_x$ is the Galois closure
of~$L_x$. Consider now an element $u$ of $(\AAf\otimes M_x)^*$ that is
equal to $p$ at one place above $p$ and equal to $1$ at all other
finite places of~$M_x$. Then $r_x(u)$, viewed as an element of
$(\AAf\otimes L_x)^*$, is equal to $p$ at exactly two places of $L_x$
above $p$ that are not in the same $\Gal(L_x/K)$-orbit, and equal to
$1$ at all other finite places of~$L_x$. It follows that $r_x(u)$ is
conjugated in $G'(\AAf)$, by some element in $G'(\Zhat)$, to the
element $g(p)$ that induces~$T_p$ (use that $G'(\ZZ_p)$ acts
transitively on the set of free rank one $\ZZ_p\otimes O_K$-submodules
of $\ZZ_p\otimes O_K^2$).  We conclude that $r_x(u)x$ is
in~$T_px$. But since $x$ is in $C_\QQ(\Qbar)$, $r_x(u)x$ is also
in~$C_\QQ(\Qbar)$. It follows that $\Gal(\Qbar/\QQ)x$ is contained in
$C_\QQ(\Qbar)\cap (T_pC_\QQ)(\Qbar)$. 
\end{proof}
Theorem~\ref{thm3.6.2} gives a lower bound for $|\Gal(\Qbar/\QQ)x|$,
whereas Theorem~\ref{thm7.2} gives an upper bound for
$|C_\QQ(\Qbar)\cap (T_pC_\QQ)(\Qbar)|$, assuming that the intersection
is finite. What we want, of course, is to show that we can choose $x$
and then $p$ such that the lower bound exceeds the upper bound, and
conclude that $C_\QQ$ and $T_pC_\QQ$ do not intersect properly. We
note that if $x$ varies over the infinite set of CM points of
$C(\Qbar)$, then $|\discr(R_x)|$ tends to infinity because there are
only finitely many orders of degree $2$ over $O_K$ with a given
discriminant, and for each such order there are only finitely many $x$
in $S(\Qbar)$ with $R_x$ isomorphic to that order.  Since our lower
bound for $|\Gal(\Qbar/\QQ)x|$ is a positive constant times a positive
power of $|\discr(R_x)|$, and our upper bound for $|C_\QQ(\Qbar)\cap
(T_pC_\QQ)(\Qbar)|$ is some fixed power of $p$, we get what we want if
we can take, for $|\discr(R_x)|$ big, $p$ of size something polynomial
in~$\log|\discr(R_x)|$. We note that
$|\discr(O_{M_x})|\leq|\discr(R_x)|^4$ because $M_x$ is the composite
of the extension $L_x$ of $K$ and its conjugate.  We also note that
the number of primes dividing $\discr(R_x)$ is at most
$\log_2(|\discr(R_x)|)$.

At this point we invoke the effective Chebotarev theorem of Lagarias,
Montgomery and Odlyzko, assuming GRH, as stated
in~\cite[Thm.~4]{Serre1} and the second remark following that
theorem. A simple computation shows that this theorem implies the
following result.
\begin{subprop}
For $M$ a finite Galois extension of $\QQ$, let $n_M$ denote its
degree, $d_M$ its absolute discriminant $|\discr(O_M)|$, and for $x$
in $\RR$, let $\pi_{M,1}(x)$ be the number of primes $p\leq x$ that
are unramified in $M$ and such that the Frobenius conjugacy class
$\Frob_p$ contains just the identity element of~$\Gal(M/\QQ)$.  Then
for $M$ a finite Galois extension of $\QQ$ and $x$ sufficiently big
(i.e., bigger than some absolute constant), and bigger than
$2(\log d_M)^2(\log(\log d_M))^2$, one has:
$$
\pi_{M,1}(x)\geq \frac{x}{3n_M\log(x)}. 
$$
\end{subprop}
This result shows that there exist infinitely many primes $p$ such that 
$C_\QQ$ and $T_pC_\QQ$ do not intersect properly. Since $C_\QQ$ is 
irreducible, it follows that, for such primes $p$, $C_\QQ$ is contained 
in~$T_pC_\QQ$. 

Assume now that $C$ is not of Hodge type. Then Proposition~\ref{prop3.5.1} 
tells us that for all primes $p$ large enough, $T_pC$ is irreducible. 
Since the correspondence $T_p$ is defined over $\QQ$, i.e., is given by a 
correspondence on $S'_\QQ$, it follows that $T_pC_\QQ$ is irreducible for 
$p$ large enough. But then we see that there exist infinitely many prime 
numbers $p$ such that $C_\QQ$ is equal to~$T_pC_\QQ$. But this is absurd, 
since by Lemma~\ref{lem3.8.1} below, 
for each $x$ in $S'(\CC)$, the Hecke orbit $\cup_{n\geq0}T_p^nx$ 
is dense in $S'(\CC)$ if $p$ is unramified in~$K$. This finishes the proof of 
Theorem~\ref{thm3.1}. 
\begin{sublem}\label{lem3.8.1}
Let $x$ be in $S'(\CC)$ and let $p$ be a prime number that is not ramified 
in~$K$. Then the Hecke orbit $\cup_{n\geq0}T_p^nx$ is dense in $S'(\CC)$ for 
the archimedean topology. 
\end{sublem}
\begin{proof}
By Lemma~\ref{lem8.4} below,
$g_0:=(\begin{smallmatrix}p&0\\0&1\end{smallmatrix})$ and $G'(\ZZ_p)$
generate $G'(\QQ_p)$ (here we use that $p$ is not ramified in~$K$).
Let now $x$ be in $S'(\CC)$, and let $(y,g)$ be a preimage of it in
$X\times G'(\AAf)$ under the quotient map for the action
by~$G(\QQ)\times G(\Zhat)$. The fact that $T_p$ is then given by right
multiplication on $G'(\AAf)$ by the element $g_0$ at the place $p$
shows that the $T_p$-orbit of $x$ is the image in $S'(\CC)$ of the
$G'(\QQ_p)$-orbit of~$(y,g)$. Let now $\Gamma$ be the subgroup of
$G'(\QQ)$ consisting of $\gamma$ such that $\gamma g$ is in~$g
G'(\Zhat)G'(\QQ_p)$. Then the $T_p$-orbit of $x$ is the image in
$S'(\CC)$ of the subset $\Gamma y$ of~$X\times\{g\}$. Now one notes
that $\Gamma$ contains a congruence subgroup of~$G'(\ZZ[1/p])$. It
follows that the intersection of $\Gamma$ with
$\SL_2(O_K[1/p])=G^\der(\ZZ[1/p])$ is dense in $G^\der(\RR)$ (for the
archimedean topology) because $G^\der$ is generated by additive
subgroups.  Since $G^\der(\RR)$ acts transitively on $X^+=(\HH^+)^2$,
the lemma is proved.
\end{proof}
\begin{sublem}\label{lem8.4}
Let $p$ be a prime that is not ramified in~$K$. Then
$g_0:=(\begin{smallmatrix}p&0\\0&1\end{smallmatrix})$ and $G'(\ZZ_p)$
generate~$G'(\QQ_p)$. 
\end{sublem}
\begin{proof}
In order to minimize notation, let $O_{K,p}$ denote $\ZZ_p\otimes
O_K$, let $K_p$ denote $\QQ_p\otimes O_K$, and let $H$ denote the
subgroup of $G'(\QQ_p)$ generated by
$g_0:=(\begin{smallmatrix}p&0\\0&1\end{smallmatrix})$ and~$G'(\ZZ_p)$.
Let $Y$ be the set of $O_{K,p}$-lattices in $K_p^2$ on which the
$K_p$-bilinear form $\psi$ given by
$(\begin{smallmatrix}0&1\\-1&0\end{smallmatrix})$ is a perfect pairing
of $O_{K,p}$-modules, up to a factor in~$\QQ_p^*$. The map
$G'(\QQ_p)\to Y$, $g\mapsto gO_{K,p}^2$ induces a bijection from
$G'(\QQ_p)/G'(\ZZ_p)$ to~$Y$. Hence, in order to prove our claim, it
suffices to show that $H$ acts transitively on~$Y$. So let $L$ be
in~$Y$.  We note that $O_{K,p}$ is either a product of two copies of
$\ZZ_p$, or the ring of integers $\ZZ_{p^2}$ in the unramified
quadratic extension $\QQ_{p^2}$ of~$\QQ_p$; in both cases, $O_{K,p}$
is a product of discrete valuation rings with uniformizer~$p$. The
theory of finitely generated modules over a discrete valuation ring
says that there exists $r$ in $\ZZ$ and $d_1$ and $d_2$ in $O_{K,p}$
such that $p^rL$ is contained in $O_{K,p}^2$ and has an
$O_{K,p}$-basis of the form $(d_1e_1,d_2e_2)$, with $(e_1,e_2)$ the
standard basis of~$O_{K,p}^2$. We note that conjugating $g_0$ by
suitable elements of $G'(\ZZ_p)$ shows that
$(\begin{smallmatrix}1&0\\0&p\end{smallmatrix})$ and
$(\begin{smallmatrix}p&0\\0&p\end{smallmatrix})$ are in $H$, and that,
in the split case, $((\begin{smallmatrix}1&0\\0&p\end{smallmatrix}),
(\begin{smallmatrix}p&0\\0&1\end{smallmatrix}))$ is in~$H$. Since the
element $d_1d_2$ of $O_{K,p}$ is the factor by which $\psi$ differs
from a perfect pairing on $p^rL$, it is actually in~$\ZZ_p$. It
follows that $(\begin{smallmatrix}d_1&0\\0&d_2\end{smallmatrix})$ is
in~$H$. This finishes the proof that $H$ is~$G'(\QQ_p)$.
\end{proof}
Let us now prove Theorem~\ref{thm3.2}. We keep the notations of the
proof Theorem~\ref{thm3.1}, and we assume again that $C$ is not of
Hodge type.  So now we may suppose moreover that $C$ contains
infinitely many CM points that have the same CM type. In particular,
we have infinitely many CM points $x$ such that $L_x$ and $\Phi_x$ are
constant, say $L$ and~$\Phi$. Of course, the orders $R_x$ are such
that $|\discr(R_x)|$ tends to infinity. The classical Chebotarev
theorem (see for example \cite[Ch.~VIII, \S4]{Lang1}) asserts that the
set of primes $p$ that are split in $M$ has natural density $1/[M:\QQ]$
(actually, Dirichlet density is good enough here). Also, recall that
the number of primes dividing some $\discr(R_x)$ is at most
$\log_2(|\discr(R_x)|)$. Hence there do exist $x$ and $p$ such that
$p$ is split in $M$, split in $R_x$, sufficiently large so that
$T_pC_\QQ$ is irreducible, and such that the lower bound for
$|\Gal(\Qbar/\QQ)x|$ of Theorem~\ref{thm3.6.2} exceeds the upper bound
for the intersection $C_\QQ(\Qbar)\cap (T_pC_\QQ)(\Qbar)$, if it is
finite. Then we have $C_\QQ=T_pC_\QQ$, hence a contradiction because
of Lemma~\ref{lem3.8.1}.

\subsection*{Acknowledgements.} It is a pleasure to thank Rutger Noot, for 
teaching a very useful course on Shimura varieties, for answering my
questions concerning them, and for his useful comments on this
manuscript.  Without the influence of Johan de Jong, Ben Moonen and
Frans Oort I would not even have started working on this subject. I
thank Paula Cohen for sending me a preliminary version
of~\cite{CohenWustholz1}, and for pointing out to me which version of
the Andr\'e-Oort conjecture was exactly needed in it. Peter
Stevenhagen gave me the reference to Stark's results
in~\cite{Stark1}. I thank Andrei Yafaev for useful remarks on this
manuscript; his numerous questions on the subject have caused me to
learn a good deal about Shimura varieties.  I thank the organizers of
the Texel conference for their excellent work. The referee deserves
much credit for pointing out a serious mistake in Section~\ref{sec3.7}
of the preprint version, and for a long list of detailed comments.
Last but not least I thank my wife Reinie Ern\'e for her influence on
this article via conversations both at work and at home, and for
letting me sleep when I am too tired to fetch a bottle for Tom at five
o'clock in the morning.

\vfill
\noindent
Bas Edixhoven\\
IRMAR\\
Campus de Beaulieu\\
35042 Rennes cedex\\
France

\newpage
\appendix
\section{Abelian surfaces with real multiplication.}
As above, $K$ is a real quadratic  field, and $O_K$ is its ring of 
integers. 
Let us describe a  bijection 
$$
S(\CC)=G(\QQ)\backslash ((\HH^{\pm})^2\times G(\AAf)/G(\Zhat)) 
\overset{\sim}{\lto} 
\{(A,\alpha)\}/\cong, 
$$
where $A$ is an abelian surface and $\alpha\colon O_K\to \End(A)$ a 
morphism of rings. By the way things have been set up, $(\HH^\pm)^2$ is 
the set of Hodge structures of type $\{(-1,0),(0,-1)\}$ on the $K$-vector 
space $K^2$, i.e., Hodge structures for which $K$ acts by endomorphisms. 

The set $G(\AAf)/G(\Zhat)$ is the set of $O_K$-lattices in~$K^2$. (By
an $O_K$-lattice in $K^2$ we mean a sub-$O_K$-module $M$ of finite
type that generates $K^2$ as a $K$-vector space). To see why the two
sets are equal, we need an adelic description of the lattices. By an
$O_K^\swedge$-lattices in $\AAf[K,]^2$ we mean a
sub-$O_K^\swedge$-module of $\AAf[K,]^2$ that is free of rank two (and
hence give the full $\AAf[K,]^2$ after tensoring with~$\QQ$). (It is
equivalent to consider sub-$O_K^\swedge$-modules of $\AAf[K,]^2$ that
are of finite type and that generate $\AAf[K,]^2$ as $K$-vector space
(or, equivalently, as $\AAf[K,]$-module).)  Now let $G(\AAf)$ act on
the set of $O_K^\swedge$-lattices in~$\AAf[K,]^2$.  This action is
transitive (use that each such a lattice is free of rank two), and the
stabilizer of the standard lattice $(O_K^\swedge)^2$ is
precisely~$G(\Zhat)$. This means that $G(\AAf)/G(\Zhat)$ is the set of
$O_K^\swedge$-lattices in~$\AAf[K,]^2$.

Let us now see why the set of $O_K$-lattices in $K^2$ is the set of
$O_K^\swedge$-lattices in~$\AAf[K,]^2$. This is always the same story,
but let me just write it down in this case (the ``classical case'' as
far as I'm concerned is for~$\ZZ$). Let $M$ be an $O_K$-lattice
in~$K^2$.  To it, we associate the $O_K^\swedge$-lattice $\Zhat\otimes
M$ in $\AAf\otimes M=\AAf[K,]^2$. In the other direction, let $N$ be a
$O_K^\swedge$-lattice in~$\AAf[K,]^2$. To $N$, we simply
associate~$N\cap K^2$.  These two maps are inverses of each other.

We can now show that $S(\CC)$ is the set of $(A,\alpha)$ up to
isomorphism. Let $(A,\alpha)$ be given. Choose an isomorphism of
$K$-vector spaces between $K^2$ and $\rH_1(A,\QQ)$. Then we get a
Hodge structure on $K^2$ and an $O_K$-lattice in $K^2$, hence an
element of $(\HH^{\pm})^2\times G(\AAf)/G(\Zhat)$, defined up to the
choice of isomorphism, i.e., up to~$G(\QQ)$. Conversely, an element of
$(\HH^{\pm})^2\times G(\AAf)/G(\Zhat)$ gives a pair $(A,\alpha)$, of
which the isomorphism class depends only on the $G(\QQ)$-orbit. So,
after all, one just has to view complex abelian varieties as given by
a $\QQ$-Hodge structure and a lattice, and use the usual stuff
regarding lattices.

One can of course do something fancy now with the category of abelian
varieties up to isogeny, and interpret $(\HH^{\pm})^2\times G(\AAf)$
as the set of isomorphism classes of $(A,\alpha,\beta,\gamma)$, with
$A$ an abelian surface up to isogeny, $\alpha$ an $K$-action on it,
$\beta$ an isomorphism of $\AAf[K,]$-modules from $\AAf[K,]^2$ to
$\rH_1(A,\AAf)$, and $\gamma$ an isomorphism of $K$-vector spaces from
$K^2$ to $\rH_1(A,\QQ)$.

\section{Polarizations.}
Why do we never have to discuss polarizability of our Hodge
structures?  Well, that's because they are in a sense only of
dimension two, just as in the case of elliptic curves. So what is in
fact true is that every complex torus of dimension two, with an action
by $O_K$, is automatically an abelian variety. Of course, this is very
standard, but I just write it down for myself, so that I understand
it, and so that I have the argument available electronically.

Consider the standard symplectic form on the $K$-vector space $K^2$: 
$$
\psi_0\colon K^2\times K^2 \lto K, \quad (x,y)\mapsto 
x^tJ y, \quad J = 
\begin{pmatrix} 0 & 1\\-1 & 0\end{pmatrix}. 
$$
Then we compose this $\psi_0$ with an arbitrary non-zero $\QQ$-linear map 
$l$ from $K$ to $\QQ$ in order to get a $\QQ$-bilinear anti-symmetric 
form: 
$$
\psi_l \colon K^2\times K^2 \overset{\psi_0}{\lto} K \overset{l}{\lto}
\QQ.
$$
In particular, one can take for $l$ the trace map; in that case, we
will denote $\psi_l$ simply by~$\psi$.

For all $g$ in $\GL_2(K)$ and all $x$ and $y$ in $K^2$, one has
$\psi_0(gx,gy)=\det(g)\psi_0(x,y)$. Hence if moreover $\det(g)$ is in
$\QQ$, one has: $\psi_l(gx,gy)=\det(g)\psi_l(x,y)$. This means that
such a $\psi_l$ is a Hodge class of weight two (or is it $-2$?) for
all $h\colon\SS_\RR\to G_\RR$ that are in $(\HH^\pm)^2$, since they
factor through the subgroup of $G$ of elements that have determinant
in~$\GG_{{\rm m},\QQ}$.

Let us now check that $\psi$ is a polarization on $(\HH^+)^2$.  So the
only condition left to check is that $x\mapsto \psi(x,h(i)x)$ should
be positive definite on $\RR\otimes K^2$. (One checks indeed that this
makes sense, in the sense that $(x,y)\mapsto\psi(x,h(i)y)$ is
symmetric.)  The fact that $\RR\otimes K=\RR^2$ means that it suffices
to check that $x\mapsto x^th(i)x$ is positive definite on $\RR^2$ for
every $h$ in~$\HH^+$.  But now note that for all $h$ in $\HH^+$ and
all non-zero $x$ in $\RR^2$, $x$ and $h(i)x$ are $\RR$-linearly
independent (interpret $h$ as giving a structure of complex vector
space on~$\RR^2$). Hence either $x\mapsto x^th(i)x$ is positive
definite, or negative definite, for all $h$ in $\HH^+$
simultaneously. So let us check just what it is for the standard $h$,
the one that sends $a+bi$ to $\bigl(\begin{smallmatrix}a & -b\\b &
a\end{smallmatrix}\bigr)$.  In that case, one has $h(i)=-J$, hence
$x^tJh(i)x=x^tx$, which is the standard inner product.

On $\HH^-$, $x^tJh(i)x$ is negative definite, since the standard $h$
there has $h(i)=J$. So we have seen that $\psi$ is a polarization on
$(\HH^+)^2$, and $-\psi$ one on $(\HH^-)^2$. On the other two
components of $(\HH^\pm)^2$ one gets polarizations by varying the map
$l$, for example by taking the composition of $l$ with multiplication
by a suitable element of $K$, i.e., an element with the right signs at
the two infinite places.

\section{Some stuff on group (schemes).}
Let $G$ denote the group scheme $\GL_2$ (over $\ZZ$, that is). Let $V$
denote its standard representation.  I use in the text that the kernel
of $G$ acting on $\Sym^2(V)\otimes\det^{-1}$ is precisely the scalar
subgroup $\Gm$ of~$G$. This can be checked as follows.  First of all,
the pairing $V\times V\to \det(V)$, $(x,y)\to xy$, is perfect.  Hence
it gives us an isomorphism between $V\otimes\det^{-1}$ and $V^*$, the
dual of~$V$. Hence we may as well consider $\End(V)=V^*\otimes V$ as
$V\otimes V\otimes\det^{-1}$. Under this isomorphism, the quotient
$\Sym_2(V)\otimes\det^{-1}$ of $V\otimes V\otimes \det^{-1}$
corresponds to the quotient of $\End(V)$ by the submodule of scalar
matrices. So we test the question there. One computes:
\begin{align*}
\begin{pmatrix}a&b\\c&d\end{pmatrix}\begin{pmatrix}1&0\\0&0\end{pmatrix} 
-\begin{pmatrix}1&0\\0&0\end{pmatrix}\begin{pmatrix}a&b\\c&d\end{pmatrix}&=
\begin{pmatrix}0&-b\\c&0\end{pmatrix}, \\
\begin{pmatrix}a&b\\c&d\end{pmatrix}\begin{pmatrix}0&0\\1&0\end{pmatrix} 
-\begin{pmatrix}0&0\\1&0\end{pmatrix}\begin{pmatrix}a&b\\c&d\end{pmatrix}&= 
\begin{pmatrix}b&0\\d-a&b\end{pmatrix}.
\end{align*} 
The condition that these two matrices are scalar give that $a=d$ and $b=c=0$. 

\section{The Shimura datum for $G^\ad$.}
For $\pi\colon G\to G^\ad$ as above, we claim that $\pi$ induces an
isomorphism (of real algebraic varieties) from $(\HH^\pm)^2$ in
$\Hom_\RR(\SS,G_\RR)$ to a conjugacy class (that we will also denote
by $(\HH^\pm)^2$) in $\Hom_\RR(\SS,G^\ad_\RR)$. Let us denote the
kernel of $\pi$ by $Z$ (it is the center of~$G$). Let $h_0$ be our
standard element in~$(\HH^\pm)^2$. Let $g$ be in $G(\RR)$, and suppose
that $\pi\circ\intaut_g\circ h_0=\pi\circ h_0$. We have to show that
$\intaut_g\circ h_0=h_0$. Here is how that goes. Define a map (i.e.,
morphism of real algebraic varieties) $z\colon \SS\to Z_\RR$ by:
$z(s)=gh_0(s)g^{-1}{\cdot}h_0(s)^{-1}$. Then, because it goes to the
center, $z$ is actually a morphism of groups. All we have to show is
that it is trivial. Well, it is trivial on the $\Gm[,\RR]$ in $\SS$,
since that one is mapped centrally in~$G_\RR$. Now the argument is
finished by noting that $Z_\RR$ is a split torus, and $\SS/\Gm[,\RR]$
is not split.

\section{Comparing various groups.}
It is not yet clear to me with which group I actually want to work.
The possibilities are: $\GL_2(K)$, $\PGL_2(K)$ and
$G'=\{g\in\GL_2(K)\;|\;\det(g)\in\QQ\}$. Just to get some idea of what
actually happens with these groups, and with the morphisms of Shimura
data between them, I think it is a good idea to make some things
explicit, such as the sets of connected components, and the finite
maps between the various Shimura varieties.

So let us first think a bit about the~$\pi_0$'s. Let's first consider 
$G$ as above. Then clearly we have: 
$$
\pi_0(S(\CC)) = \GL_2(K)^+\backslash \GL_2(\AAf[K,])/\GL_2(O_K^\swedge). 
$$
But this set is the set of isomorphism classes of locally free rank 2
$O_K$-modules $M$ with an orientation on $\det(M)=\Lambda^2_{O_K}M$ at
the two infinite places. Now each locally free rank two module over
$O_K$ is isomorphic to one of the form $O_K\oplus L$ (show first that
it is decomposable by choosing a one dimensional $K$-sub-vector space
in $\QQ\otimes M$; then show that $M$ has a nowhere vanishing
element).  Of course, $L$ is determined by $M$ since one has
$\det(M)=L$.  It follows that:
$$
\text{$\pi_0(S(\CC)) = \Pic(O_K)^+$, the strict class group of $K$.}
$$
Let us now consider $\pi_0(S^\ad(\CC))$. We have: 
$$
\pi_0(S^\ad(\CC)) = \PGL_2(K)^+\backslash
\PGL_2(\AAf[K,])/\PGL_2(O_K^\swedge).
$$
This we recognize as the set of isomorphism classes of $\PP^1$-bundles
on $S:=\Spec(O_K)$, locally trivial in the Zariski topology, with an
orientation at the two infinite places (it does not seem a complete
tautology, the correspondence with the Zariski $\PP^1$-bundles,
namely, it says more directly something as: trivial over $K$, and over
every completion).  Anyway, let us show that each $\PP^1$-bundle $X$
on $S$ comes from a locally free rank two bundle on~$S$. Just note
that each element of $X(K)$ extends to one in~$X(S)$. An element in
$X(S)$ gives an invertible $\calO_X$-module $\calL$ that has degree
one on each fibre, hence with $p_*\calL$ a rank two bundle on
$S$. Then one checks that $X$ is isomorphic, over $S$, to
$\PP(p_*\calL)$ (it is easy to see that $X$ is the Grassmannian of
locally free rank one quotients of $p_*\calL$).  This has an
interpretation in the long exact sequence coming from the short exact
sequence of Zariski sheaves on $S$:
$$
1\lto\Gm[,S]\lto\GL_{2,S}\lto\PGL_{2,S}\lto 1. 
$$
What is quite nice in this situation is that $S$ is of dimension one,
hence that $\rH^2(S,\Gm)=0$, which explains the observation above.
Now that we know that each element of $\pi_0(S^\ad(\CC))$ comes from a
locally free rank two $O_K$-module, we want to know when two such
modules give isomorphic $\PP^1$-bundles. Well, the Grassmannian
interpretation says that that happens if and only if the two modules
are isomorphic up to twist by an invertible $O_K$-module. Hence:
$$
\FF_2\otimes\Pic(O_K)=\{\text{$\PP^1$-bundles on $\Spec(O_K)$}\}_{/\cong}. 
$$
Let us now consider orientations. Note that $\Aut_{O_K}(O_K\oplus L)$
maps surjectively, under $\det$, to $\Aut_{O_K}(L)=O_K^*$, and that
doubles in $\Pic(O_K)$ have a canonical orientation. It follows that:
$$
\pi_0(S^\ad(\CC)) = \Pic(O_K)^+/2\Pic(O_K), 
$$
and that: 
$$
\#\pi_0(S^\ad(\CC)) = 
\begin{cases} 
\#\FF_2\otimes\Pic(O_K) &  \text{if $N(O_K^*)=\{1,-1\}$,}\\
2\#\FF_2\otimes\Pic(O_K) & \text{if $N(O_K^*)=\{1\}$.}
\end{cases}
$$
Let us now say something about the map $S(\CC)\to S^\ad(\CC)$.  This
makes it necessary to know things about the morphism of group schemes
$\GL_{2,O_K}\to\PGL_{2,O_K}$.  We would like to know that this
morphism is surjective for the Zariski topology. For this, it suffices
to show that on a scheme $S$, an $S$-automorphism of $\PP^1_S$ is
induced, locally on $S$, by an element of~$\GL_2(S)$. Now use that for
any scheme $T$, to give an element of $\PP^1(T)$ is to give an
invertible $O_T$-module with two sections that generate it. Let $g$ be
an $S$-automorphism of $\PP^1_S$.  Then $g^*\calO(1)$ is of the form
$p^*\calL\otimes\calO(1)$ for some invertible
$\calO_S$-module~$\calL$. Since $\calL$ is locally trivial, we get
what we want.

Hence: the morphism $\GL_2(O_K^\swedge)\to\PGL_2(O_K^\swedge)$ is
surjective (use that $O_K^\swedge$ is the product of the completions
at all finite places and that one has the surjectivity for each such
completion).  And: $\GL_2(\AAf[K,])\to\PGL_2(\AAf[K,])$ is surjective
(just use what elements of $\AAf[K,]$ look like, or use that to give a
point of a scheme with values in $\AAf$ is to give, for each $p$, a
point with values in $\QQ_p$, such that for almost all $p$ the point
comes from a point with values in~$\ZZ_p$). The more difficult thing
that remains now is the question of surjectivity of the morphism
$\GL_2(O_K)\to\PGL_2(O_K)$.

It follows that $S(\CC)\to S^\ad(\CC)$ is surjective. The stabilizer
of $(\HH^+)^2\times\{1\}$ in $\GL_2(K)$ and $\PGL_2(K)$ are
$\GL_2(O_K)^+$ and $\PGL_2(O_K)^+$, respectively. So let us find out
what the cokernel of $\GL_2(O_K)^+\to \PGL_2(O_K)^+$ is.

An automorphism of $\PP:=\PP^1_S$ is given by an invertible
$\calO_\PP$-module of degree one together with two generating
sections. Such a module is of the form $p^*\calL\otimes\calO(1)$. But
then we have the condition that
$p_*p^*\calL\otimes\calO(1)=\calL\oplus\calL$ is generated by two
global sections. This can be done if and only if
$\calL\oplus\calL\cong\calO\oplus\calO$, i.e., if and only if
$\calL^{\otimes2}\cong\calO$. This explains that we have an exact
sequence:
$$
1\lto O_K^*\lto\GL_2(O_K)\lto\PGL_2(O_K)\lto\Pic(O_K)[2]\lto 0
$$
Likewise, one gets: 
$$
1\lto O_K^*\lto\GL_2(O_K)^+\lto\PGL_2(O_K)^+\lto\Pic(O_K)[2]\lto 0. 
$$
We conclude that $S(\CC)^0\to S^\ad(\CC)^0$ is the quotient for a
faithful action by the group $\Pic(O_K)[2]$, where $S(\CC)^0$ and
$S^\ad(\CC)^0$ are the standard irreducible components of $S(\CC)$
and~$S^\ad(\CC)$.  One computes directly that the map $S(\CC)\to
S^\ad(\CC)$ is the quotient for a faithful action by the group
$K^*\backslash\AA_K^*/O_K^{\swedge,*}$, i.e., by~$\Pic(O_K)$.

Let us now do some comparing between $S$ and $S'$, with $S'$ coming
from the Shimura datum with the group~$G'$. We first remark that
$G'(\RR)$ is the set of $(g_1,g_2)$ in $\GL_2(\RR)^2$ such that
$\det(g_1)=\det(g_2)$.  This means that the $G'(\RR)$ conjugacy class
of morphisms from $\SS$ to $G'_\RR$ that we deal with is:
$$
X':=(\HH^+)^2\coprod(\HH^-)^2=(\HH^2)^\pm. 
$$
Hence we have: 
$$
S'(\CC) = G'(\QQ)\backslash(X'\times G'(\AAf)/G'(\Zhat)),\quad 
\pi_0(S'(\CC)) = G'(\QQ)^+\backslash G'(\AAf)/G'(\Zhat). 
$$
This last set is the set of isomorphism classes of triplets
$(M,\phi,\alpha)$ with $M$ a locally free $O_K$-module of rank two,
$\phi\colon O_K\to\det(M)$ an isomorphism, and $\alpha$ an orientation
on $\RR\otimes_\ZZ\det(M)$ that induces plus or minus the standard
orientation on $\RR\otimes O_K=\RR\times\RR$ via~$\phi$.  Since every
such triplet is isomorphic to $(O_K^2,\id,(+,+))$, we see:
$$
\text{$S'(\CC)$ is connected.}
$$
What about the map $S'(\CC)\to S(\CC)$?  It suffices to look at what
happens on $(\HH^2)^+\times\{1\}$. The stabilizer of this in $G'(\QQ)$
is simply $\SL_2(O_K)$, and the stabilizer in $G(\QQ)$ is
$\GL_2(O_K)^+$. Hence the map $S'(\CC)\to S(\CC)^+$ is the quotient
for the faithful action by the group $O_K^{*,+}/O_K^{*,2}$, i.e.,
totally positive global units modulo squares of global units.

\section{Some stuff on bilinear forms and field extensions.}
Let $k\to K$ be a finite field extension, say of degree~$d$. Let $V$
be a finite dimensional $K$-vector space, say of dimension~$n$.  Let
$X$ denote the $k$-vector space of $k$-bilinear forms $b\colon V\times
V\to k$ such that $b(ax,y)=b(x,ay)$ for all $x$ and $y$ in $V$ and all
$a$ in~$K$. (I.e., the maps $x\mapsto ax$ are required to be
self-adjoint.) We want to relate $X$ to the set $Y$ of $K$-bilinear
forms on~$V$.

Let $l\colon K\to k$ be a surjective $k$-linear map (for example, one can 
take the trace map if $k\to K$ is separable). Then we have a map: 
$$
L\colon Y\lto X, \quad b\mapsto l\circ b. 
$$
Indeed, for $b$ in $Y$ we have: $(l\circ b)(ax,y)=l(b(x,ay))=(l\circ
b)(x,ay)$. The map $L$ is injective, since, for $b$ a $K$-bilinear
form, the image of $b$ is either $0$ or~$K$.

Let us now assume that $k\to K$ is separable. Then one computes that
both $X$ and $Y$ are of dimension $n^2d$ over $k$ (of course, for $Y$
this is clearly true without the separability assumption; for $X$, one
uses this assumption in order to reduce to the case $K=k^d$ via base
change from $k$ to some algebraic closure for example). Hence we
conclude that the map $L$ above is bijective. (I did not bother to
check if this is still true without the separability.)  So we have the
following result.
\begin{prop}
Let $k$ be a field, and $K$ a finite separable $k$-algebra. Let
$l\colon K\to k$ be a surjective $k$-linear map (for example the trace
map).  Let $V$ be a finitely generated projective $K$-module. Then for
every $k$-bilinear form $b\colon V\times V\to k$ such that
$b(ax,y)=b(x,ay)$ for all $x$ and $y$ in $V$ and all $a$ in $K$ there
exists a unique $K$-bilinear $b'\colon V\times V\to K$ such that
$b=l\circ b'$. With this notation, $b$ is symmetric (antisymmetric) if
and only if $b'$ is so.
\end{prop}
\begin{proof}
It only remains to prove that $b$ is symmetric (antisymmetric) if and
only if $b'$ is so. For $b$ as above, let $b^t$ denote its adjoint,
i.e., $b^t(x,y)=b(y,x)$; we will use the same notation for elements
of~$Y$.  Then one has $(b^t)'=(b')^t$. Now $b$ is symmetric if and
only if $b^t=b$, and $b$ is antisymmetric if and only if
$b^t=-b$. Hence the result.
\end{proof}
The next result gives a construction of the inverse of $L$, if one
takes $l$ to be the trace map.
\begin{prop}
Let $k$ be a field, and $K$ a finite separable $k$-algebra. Let $V$ be
a finitely generated projective $K$-module, and $b\colon V\times V\to
k$ a $k$-bilinear map such that $b(ax,y)=b(x,ay)$ for all $x$ and $y$
in $V$ and all $a$ in $K$. Because of the separability, we have a
natural isomorphism of $K$-algebras: $K\otimes_kK=K\times K'$, where
we view $K\otimes_kK$ as a $K$-algebra via the first factor. This
decomposition gives a decomposition of $K$-modules:
$K\otimes_kV=K\otimes_kK\otimes_KV=V\oplus V'$ with
$V'=K'\otimes_KV$. Let $b_K$ denote the $K$-bilinear form on
$K\otimes_kV$ obtained by extension of scalars. Then the decomposition
of $K\otimes_kV$ in $V$ and $V'$ is orthogonal for $b_K$, and $b'$ is
the restriction to $V$ of $B_K$. In particular, one has $b=\trace\circ
b'$.
\end{prop}
Let us now note the special case where $V$ is of dimension two. In
that case, the $K$-vector space $Y$ of antisymmetric $K$-bilinear
forms is of dimension one, hence one gets the following corollary,
which is of interest for Hilbert modular varieties.

\begin{cor}
Let $k$ be a field, and $k\to K$ a finite separable $k$-algebra. Let
$V$ be a free $K$-module of rank two. Let $\psi_0\colon V\times V\to
K$ be a non-degenerate alternating $K$-bilinear form. Then for every
alternating $k$-bilinear form $\psi\colon V\times V\to k$ such that
$\psi(ax,y)=\psi(x,ay)$ for all $x$ and $y$ in $V$ and all $a$ in $K$,
there exists a unique $b$ in $K$ such that
$\psi(x,y)=\trace(b\psi_0(x,y))$ for all $x$ and $y$ in~$V$.
\end{cor}

\section{Moduli interpretation for the symplectic group.}
For details, see \cite[Sections 1, 4]{Deligne2}.  Just in this
section, let $G$ denote the group of symplectic similitudes of rank
$2n$. More precisely, let $n\geq 0$ be an integer, and let $G$ denote
the group of automorphisms of the $\ZZ$-module $\ZZ^{2n}$ that
preserve, up to scalar multiple, the standard symplectic form, i.e.,
the form given by the matrix
$(\begin{smallmatrix}0&1\\-1&0\end{smallmatrix})$.  Let $X:=\HH_n^\pm$
the set of $h\colon\SS\to G_\RR$ that are Hodge structures of weight
$-1$ such that $\psi$ is a polarization up to a sign.  Then this $X$
is one $G(\RR)$-conjugacy class and it is called the Siegel double
space. Let us consider:
$$
A_n(\CC):=G(\QQ)\backslash(X\times G(\AAf)/G(\Zhat)). 
$$
What we want to show is that $A_n(\CC)$ is the set of isomorphisms
classes of pairs $(A,\lambda)$ of principally polarized abelian
varieties of dimension~$n$. We already know what the interpretation of
$X$ is: it is the set of Hodge structures of weight $-1$ such that
$\psi$ is a polarization up to a sign. Let us now interpret
$G(\AAf)/G(\Zhat)$. Consider the action of $G(\AAf)$ on the set of
lattices in $\AAf^2$. The stabilizer of the standard lattice $\Zhat^2$
is $G(\Zhat)$. Hence $G(\AAf)/G(\Zhat)$ is the set of lattices of the
form $x\Zhat^2$, with $x$ in~$G(\AAf)$. We claim that this is the set
of lattices $L$ on which a suitable multiple of $\psi$ induces a
perfect pairing. For $x$ in $G(\AAf)$ we have:
$\psi(xu,xv)=\mu(x)\psi(u,v)$, which proves that $\mu(x)^{-1}\psi$ is
a perfect pairing on~$x\Zhat^2$. On the other hand, let $L$ be a
lattice and $a$ in $\AAf^*$ be such that $a\psi$ is a perfect pairing
on~$L$.  Then take a $\Zhat$-basis $l_1,\ldots,l_{2n}$ of $L$ such
that $a\psi$ is in standard form, i.e., given by the matrix
$(\begin{smallmatrix}0&1\\-1&0\end{smallmatrix})$.  Then the element
$x$ of $\GL_{2n}(\AAf)$ with $xe_i=l_i$ is in $G(\AAf)$. This finishes
the proof of the fact that $G(\AAf)/G(\Zhat)$ is the set of lattices
on which a multiple of $\psi$ is perfect.

Let us now describe the constructions that give $A_n(\CC)$ the
interpretation as the set of isomorphism classes of abelian varieties
of dimension $n$, with a principal polarization.

Suppose $(A,\lambda)$ is given. Then choose an isomorphism
$f\colon\QQ^{2n}\to\rH_1(A,\QQ)$ such that $\psi$ corresponds to a
multiple of~$\lambda$ (such an $f$ is unique up to an element
of~$G(\QQ)$). Let $x$ be the element of $X$ that is given by the Hodge
structure on $\QQ^{2n}$ induced from $A$ via~$f$. Let $L$ in
$G(\AAf)/G(\Zhat)$ be the lattice corresponding to $\ZZ^{2n}$ via~$f$.
The class of $(x,L)$ modulo $G(\QQ)$ depends only on the isomorphism
class of $(A,\lambda)$.

Suppose now that we have $(x,L)$ in $X\times G(\AAf)/G(\Zhat)$.  Then
let $A$ be $(\RR\otimes L)/L$ with the complex structure given by the
Hodge structure corresponding to~$x$. Let $a$ be the element of
$\QQ^*$ such that $a\psi$ is perfect on $L$ (this fixes $a$ up to
sign) and is a polarization $\lambda$ on $A$ (this fixes the
sign). For $g$ in $G(\QQ)$, multiplication by $g$ gives an isomorphism
from $(A,\lambda)$ to the $(A',\lambda')$ obtained from $(gx,gL)$.

Let us end with a remark which is just a reminder to myself. 
\begin{rem}
Let $V$ be a free finitely generated $\ZZ$-module, with
$h\colon\SS\to\GL(V)_\RR$ a Hodge structure of type
$(-1,0),(0,-1)$. Let $A:=(\RR\otimes V)/V$ be the associated complex
torus. Then the dual complex torus corresponds to the Hodge structure
$t\mapsto(h(t)^\vee)^{-1}N(t)=h(\bar{v})^\vee$ on~$V^\vee$. In other
words, the dual of $A$ is $(\RR\otimes V^\vee)/V^\vee$, with the
complex structure on $\RR\otimes V^\vee=(\RR\otimes V)^{\vee_\RR}$
such that $z$ in $\CC$ acts as~$\bar{z}^\vee$. In order to prove this,
one notes that the tangent space of $A^t$ is $\rH^1(A,\calO_A)$, which
is naturally $\CC$-anti-linearly isomorphic to $\rH^0(A,\Omega^1_A)$,
which is the dual of the tangent space of $A$ at zero.
\end{rem}

\section{Moduli interpretation of $S'(\CC)$.}
Let us recall: 
$$
S'(\CC) = G'(\QQ)\backslash(X'\times G'(\AAf)/G'(\Zhat)). 
$$
\begin{prop}
The Shimura variety $S'_\QQ$ is the moduli space of triplets
$(A,\alpha,\lambda)$with $A$ an abelian surface, $\alpha\colon O_K\to
\End(A)$ a ring morphism, and $\lambda\colon A\to A^*$ a principal
$O_K$-polariza\-tion.
\end{prop}
First of all, we have to explain what $A^*$ is, and what we call a
principal $O_K$-polarization. Let us begin with $A^*$: it is the dual
of $A$ in the category of abelian varieties with $O_K$-action. More
precisely, since for $A$ an abelian variety the dual is defined to be
$A^t:=\Ext^1(A,\Gm)$, we put:
$$
A^* := \Ext^1_{O_K}(A,O_K\otimes \Gm) = \delta\otimes_{O_K}\Ext^1(A,\Gm) 
= \delta\otimes_{O_K}A^t, 
$$
with 
$$
\delta:=\Hom_{O_K}(\Hom_\ZZ(O_K,\ZZ),O_K)
$$
the different of the extension $\ZZ\to O_K$. In order to prove the
above equalities, it is useful to note that for $A\to B$ a morphism of
rings, for $M$ a $B$-module and $N$ and $A$-module, one has the
adjunction:
$$
\Hom_B(M,\Hom_A(B,N)) = \Hom_A({}_AM,N), 
$$
where ${}_AM$ denotes the $A$-module given by~$M$. Then one uses that
for $B$ locally free of finite rank as $A$-module one has
$\Hom_A(B,N)=B^{\vee_A}\otimes_AN$, with $B^{\vee_A}=\Hom_A(B,A)$ the
$A$-dual of~$B$. And then one uses that for $P$ a finitely generated
locally free $B$-module one has:
$$
\Hom_B(M,P\otimes_AN) = P\otimes_B\Hom_B(M,B\otimes_AN). 
$$
This establishes: 
$$
\Hom_B(M,B\otimes_AN) = \delta\otimes_B\Hom_A(M,N), 
$$
with $\delta=(B^{\vee_A})^{\vee_B}$. Deriving with respect to $N$ then 
gives: 
$$
\Ext^i_B(M,B\otimes_AN) = \delta\otimes_B\Ext^i_A(M,N). 
$$
This explains the reason that $A^*$ occurs in this context. In our
context, $A=\ZZ$ and $B=O_K$, so that $\delta$ is an ideal in $O_K$,
since the trace map $\trace\colon O_K\to\ZZ$ gives an injective
morphism $O_K\to(O_K)^{\vee_\ZZ}$ (in fact, the trace map from $B$ to
$A$ always gives a morphism $\delta\to B$, but it might be zero). This
gives an isogeny:
$$
A^* = \delta\otimes_{O_K}A^t\lto A^t. 
$$
We define an $O_K$-polarization to be an $O_K$-morphism $\lambda\colon
A\to A^*$ such that the induced morphism from $A$ to $A^t$ is a
polarization; $\lambda$ is called principal if it is an isomorphism
(note that this means that the induced morphism $A\to A^t$ is not an
isomorphism, since $O_K$ is ramified over $\ZZ$ (we do suppose that
$K$ is a field, after all).

Let us now turn to the proof of the proposition above that gives the
moduli interpretation of~$S'_\QQ$. So we want to show that $S'(\CC)$
is the set of isomorphism classes of triplets $(A,\alpha,\lambda)$
over~$\CC$. In Hodge theoretical terms, such triplets are given by
triplets $(V,h,\psi)$ with $V$ a locally free $O_K$-module of rank
two, $h\colon\SS\to(\GL_\ZZ(V))_\RR$ a Hodge structure of type
$(-1,0),(0,-1)$, and $\psi\colon V\times V\to O_K$ a perfect
antisymmetric $O_K$-bilinear form such that $\trace\circ\psi\colon
V\times V\to\ZZ$ is a polarization.  Note that for such a triplet
$(V,h,\psi)$, the pair $(V,\psi)$ is isomorphic to the standard pair
$(O_K\oplus O_K,(\begin{smallmatrix}0&1\\-1&0\end{smallmatrix}))$.
The proof of the proposition can now be easily described.  As in the
last section, one shows that $G'(\AAf)/G'(\Zhat)$ is the
$O_K$-lattices in $K^2$ on which
$\psi=(\begin{smallmatrix}0&1\\-1&0\end{smallmatrix})$ induces a
perfect pairing of $O_K$-modules, up to a factor in~$\QQ^*$. The space
$X'$ is the set of Hodge structures. As in the last section, one shows
that $G'(\AAf)/G'(\Zhat)$ is the $O_K$-lattices in $K^2$ on which
$\psi=(\begin{smallmatrix}0&1\\-1&0\end{smallmatrix})$ induces a
perfect pairing of $O_K$-modules, up to a factor in~$\AAf^*$.  The
space $X'$ is the set of Hodge structures of type $(-1,0),(0,-1)$ on
the $K$-vector space $K^2$ such that, up to sign, $\trace\circ\psi$ is
a polarization. f type $(-1,0),(0,-1)$ on the $K$-vector space $K^2$
such that, up to sign, $\trace\circ\psi$ is a polarization. After
these remarks one simply follows the lines of the proof above of the
modular interpretation for the symplectic group.  Anyway, for details,
one can consult \cite[4.11]{Deligne2}.

Let us end by stating that the multiplier character $\mu\colon
G'\to\Gm$ is the determinant (view $G'$ as a subgroup of
$\Res_{O_K/\ZZ}\GL_{2,O_K}$ and $\Gm$ as a subgroup of
$\Res_{O_K/\ZZ}\Gm[,O_K]$. More precisely, for all $g$ in $G'(\QQ)$
and all $x$ and $y$ in $K^2$ we have
$(\trace\circ\psi)(gx,gy)=\det(g)(\trace\circ\psi)(x,y)$.

\section{A remark on Mumford-Tate groups.}
What I want to say is that to an isomorphism class of $\QQ$-Hodge
structures one can associate its Mumford-Tate group. Namely, if $V$
and $V'$ are isomorphic $\QQ$-Hodge structures, and if $f$ and $f'$
are isomorphisms from $V$ to $V'$, then $f'=fg$ with $g$ an
automorphism of $V$. But then $g$ centralizes the Mumford-Tate group
in $\GL(V)$. Hence $f$ and $f'$ induce the same isomorphism from
$\MT(V)$ to~$\MT(V')$.  For example, the functor $V\mapsto \Aut(V)$
does not have this property.

The same argument shows that a point $P$ on a Shimura variety 
$\Sh_K(G,X)(\CC)$ defines an algebraic group $\MT(P)$, with a given 
$G(\QQ)$-conjugacy class of embeddings in~$G$. 

\section{On computing the generic Mumford-Tate group on $S'(\CC)$.}
First note that for all $h=(h_1,h_2)\colon \CC^*\to\GL_2(\RR)^2$ in
$X'$ one has $\det(h_1(z))=\det(h_2(z))$ for all~$z$. This shows that
$\MT$ is contained in~$G'_\QQ$.  The locally constant sheaf $V$
becomes constant on~$X'$. Hence $\MT_\RR$ contains all
$h(\CC^*)\subset\GL_2(\RR)^2$ for the $h$ in~$X'$.  In particular, it
contains all conjugates under $G'(\RR)$ of those images.  but then it
contains all $(x,x)$, all $(yxy^{-1},x)$, hence all
$(xyx^{-1}y^{-1},1)$, etc. It follows that $\MT=G'_\QQ$.

\section{Other remarks on Mumford-Tate groups.}
We have defined the Mumford-Tate group $\MT(V)$ of a $\QQ$-Hodge
structure $V$ given by a morphism $h\colon \SS\to\GL(V)_\RR$ to be the
smallest algebraic subgroup $H$ of $\GL(V)$ such that $h$ factors
through~$H_\RR$. This is not the usual definition, perhaps. The usual
definition is to take $\MT'(V)$, the smallest subgroup $H$ of
$\GL(V)\times\Gm$ such that $h'\colon\SS\to\GL(V)_\RR\times\Gm[,\RR]$
factors through $H_\RR$, where $\SS\to\Gm[,\RR]$ corresponds
to~$\QQ(1)$. The difference between the two choices is that $\MT'(V)$
keeps track of weights, whereas $\MT(V)$ doesn't.  The Tannakian
interpretation of $\MT(V)$ is that it is the automorphism functor of
the fibre functor ``forget Hodge structure'' on the tensor category
generated by~$V$. For $\MT'(V)$, one considers the tensor category
generated by $V$ and~$\QQ(1)$. Yet another (of course related)
characterization is that $\MT'(V)$ seems to be the biggest subgroup of
$\GL(V)\times\Gm$ that fixes all elements of type $(0,0)$ in
$\QQ$-Hodge structures of the form $V^{\otimes n}\otimes(V^*)^{\otimes
m}\otimes\QQ(p)$.  For this, see Deligne-Milne-Ogus-Shih. In the same
way, $\MT(V)$ is characterized by the fact that it stabilizes all
lines generated by Hodge classes (i.e., classes of some type $(p,p)$)
in $\QQ$-Hodge structures of the form $\oplus_iV^{\otimes
n_i}\otimes(V^*)^{\otimes m_i}$.

Since I did not find this explicitly written (but I haven't looked
very much, I should say), let me write a proof. So let $H$ be the
intersection of the stabilizers of such lines. Let us first prove that
$\MT(V)\subset H$.  So let $t$ in some $T=\oplus_iV^{\otimes
n_i}\otimes(V^*)^{\otimes m_i}$ be of some type $(p,p)$. Then $\RR
t\subset T_\RR$ is fixed by $\SS$, hence by~$\MT(V)$. This proves that
$\MT(V)\subset H$. Let's now prove that $H\subset\MT(V)$. Now we use
Chevalley's result: every subgroup of $\GL(V)$ is the stabilizer of a
line in some finite dimensional representation of $\GL(V)$, plus the
fact that each finite dimensional representation of $\GL(V)$ is
contained in a representation of the form $\oplus_iV^{\otimes
n_i}\otimes(V^*)^{\otimes m_i}$ (I will even give proofs for these two
facts below, since I do not like the proof given in DMOS). Anyway, let
$t$ in some $T$ be such that $\MT(V)$ is the stabilizer of $\QQ
t$. Then $\RR t$ is fixed by $\SS$, hence $t$ is of some type $(p,p)$.
(Use for example that the norm $\SS\to\Gm[,\RR]$ generates
$\Hom(\SS,\Gm[,\RR])$.)

As I said, I do not like the proof of parts (a) and (b) of
Proposition~3.1 in Chapter~I of DMOS. So I give one.

\begin{thm}
Let $G$ be an affine algebraic group over a field $k\supset\QQ$. Let
$H$ be an algebraic subgroup of $G$, and $V$ a finite dimensional
faithful representation of~$G$. Then there exists a line $L$ in some
representation of $G$ of the form $\oplus_iV^{\otimes
n_i}\otimes(V^*)^{\otimes m_i}$, such that $H$ is the stabilizer
of~$L$.
\end{thm}
\begin{proof}
First of all, we may and do suppose that $G=\GL(V)$. The idea is now
the following: let $G$ act on itself by right translation; then $G$
acts on $k[G]$, and $H$ is the stabilizer of the ideal $I_H$; then use
that $I_H$ is finitely generated, and that $k[G]$ is locally
finite. Let us first write down what $k[G]$ is, as a $G$-module via
right translation on~$G$.  Well,
$$
k[G]=k[\End(V)][1/\det]=\Sym_k(\End(V)^*)[1/\det]=\Sym_k(V^d)[1/\det], 
$$
where the last equality comes from the fact that $\End(V)^*$, as
$G$-module given by right translation on $\End(V)$, is simply $V^d$,
where $d$ is of course the dimension of $V$ (note that the $G$-action
on $\End(V)^*$ extends to an $\End(V)$-action). Also, note that $\det$
is in $\Sym^d(\End(V)^*)$, and that we have
$g{\cdot}\det=\det(g)\det$. The scalar subgroup $\Gm$ of $G$ induces a
$\ZZ$-grading on~$k[G]$. We have:
$$
k[G]_i=\bigcup_j k[\End(V)]_{i+dj}{\det}^{-j}, 
$$
and: 
\begin{align*}
&k[\End(V)]_i = \Sym^i(\End(V)^*) = \Sym^i(V^d)\subset(V^d)^{\otimes i} 
= (V^{\otimes i})^{d^i}, \\
&k{\cdot}\det = \Lambda^dV\subset V^{\otimes d}, \\
&k{\cdot}{\det}^{-1} = (\Lambda^dV)^* \subset(V^*)^{\otimes d}. 
\end{align*}
This describes $k[G]$ as $G$-module. Let $f_1,\ldots,f_r$ be a finite
set of generators of the ideal $I_H$ of~$k[G]$. Let $W\subset k[G]$ be
a finite dimensional sub-$G$-module containing the~$f_i$. Then $H$ is
the stabilizer of the subspace $I_H\cap W$ of $W$, hence of the line
$\Lambda^n(I_H\cap W)\subset\Lambda^n(W)$, with $n=\dim(I_H\cap
W)$. Now note that $W$ is a subrepresentation of a representation of
the form $\oplus_iV^{\otimes n_i}\otimes(V^*)^{\otimes m_i}$.
\end{proof}
\begin{rem}
If we allow subquotients of the 
$\oplus_iV^{\otimes n_i}\otimes(V^*)^{\otimes m_i}$, then we can drop the 
hypothesis that $k$ is of characteristic zero. 
\end{rem}
\begin{rem}
If $H$ contains the scalars in $G=\GL(V)$, then one can take $L$ to be
in some representation of the form $(V^{\otimes n})^m$. To prove this,
consider the Zariski closure $\ol{H}$ of $H$ in $\End(V)$, and use
that it is a cone.
\end{rem}
Just for fun, let us look at some examples in $G:=\GL_2$. 
The Borel subgroup 
$B:=\{(\begin{smallmatrix}* & * \\0 & *\end{smallmatrix})\}$ is the stabilizer 
of the line generated by $(1,0)$ in $V:=k^2$. The subgroup 
$\{(\begin{smallmatrix}1 & * \\0 & *\end{smallmatrix})\}$ is the stabilizer 
of $k(1,(1,0))$ in $k\oplus V$. The subgroup 
$\{(\begin{smallmatrix}* & * \\0 & 1\end{smallmatrix})\}$ is the stabilizer 
of $k(1,(0,1)^*)$ in $k\oplus V^*$. The subgroup 
$\{(\begin{smallmatrix}t & 0 \\0 & t\end{smallmatrix})\}$ is the stabilizer 
of $k((1,0),(0,1))$ in $V\oplus V$. The subgroup 
$\{(\begin{smallmatrix}* & 0 \\0 & *\end{smallmatrix})\}$ is the stabilizer 
of the two-dimensional subspace of the $((x,0),(0,y))$ in $V\oplus V$; note 
that the proof above gives the same result. Finally, the trivial subgroup 
$\{(\begin{smallmatrix}1 & 0 \\0 & 1\end{smallmatrix})\}$ is the stabilizer 
of $(1,(1,0),(0,1))$ in $k\oplus V\oplus V$. 

\section{Modular interpretation of $T_p$.}
Let $A$ be a complex abelian surface with multiplication by $O_K$ and
with a principal $O_K$-polariza\-tion $\lambda\colon A \to A^*$. Let
$H$ be an $O_K$-submodule of $A[p](\CC)$ that is free of rank
one. Then we claim that $p\lambda$ induces a principal
$O_K$-polarization on~$A/H$. So how does this work? Write
$p\lambda=\pi_2\pi_1$, with $\pi_1\colon A\to B$ the quotient
by~$H$. Then $(A/H)^*$ is the quotient of $A$ by $\ker(\pi_2^*)$.  So
we have to see that $\ker(\pi_2^*)=\ker(\pi_1)$. Since both have the
same number of elements, it suffices to see that one is contained in
the other.  Since $\ker(\pi_1)$ is maximal isotropic for the pairing
$e_{\lambda,p}$ that $\lambda$ induces on $A[p]$, it suffices to see
that $\ker(\pi_2^*)$ and $\ker(\pi_1)$ are orthogonal for that
pairing. That results from standard things about such pairings coming
from expressions like $A^*=\Ext^1(A,O_K\otimes\Gm)$.

The general statement is this: let $f\colon A\to B$ and $g\colon B\to
C$ be isogenies of abelian varieties with multiplications by~$O_K$.
Let $h:=gf$. Then we have a short exact sequence:
$$
0\lto\ker(f)\lto\ker(h)\lto\ker(g)\lto 0. 
$$
Applying $\Hom(\cdot,O_K\otimes\Gm)$ gives an isomorphism of short 
exact sequences: 
$$
\begin{CD}
0 @>>> \ker(g^*) @>>> \ker(h^*) @>>> \ker(f^*) @>>> 0 \\
@VVV @VVV @VVV @VVV  @VVV \\
0 @>>> \ker(g)^* @>>> \ker(h)^* @>>> \ker(f)^* @>>> 0. 
\end{CD}
$$
The fact that the map from $\ker(g^*)$ to $\ker(f)^*$ is zero means
that $\ker(f)$ and $\ker(g^*)$ are orthogonal for the paring induced
by $h$ between $\ker(h)$ and~$\ker(h^*)$.

\section{Some stuff on orders in finite separable $\QQ$-algebras.} 
I need lower bounds for orders of Picard groups of certain orders in
certain CM fields. Therefore, some general theory should be quite
useful.

Let $\QQ\to K$ be a finite separable $\QQ$-algebra. Then $K$ is a
finite product of number fields, say $K=K_1\times\cdots\times K_m$,
and the integral closure of $\ZZ$ in $K$ is then the product of the
maximal orders of the~$K_i$. Let $R\subset K$ be an order in $K$,
i.e., a subring of $K$ with $\QQ\otimes R=K$ and which is finitely
generated as a $\ZZ$-module.  Then $R$ is contained in $O_K$ since the
elements of $R$ are integral over~$\ZZ$, and $O_K/R$ is a finite
additive group, since $R$ and $O_K$ are free $\ZZ$-modules of the same
finite rank. Consider ideals of $R$ that are also
$O_K$-ideals. Clearly a lot of such ideals do exist: for every $n$ in
$\ZZ$ that annihilates $O_K/R$, we have the example~$nO_K$. The sum of
a family of such ideals is again one such, hence there exists a unique
maximal such ideal, called the conductor of $R$ (relative to~$O_K$).
I don't think that we will use this conductor so much, since we want
estimates in terms of the discriminant of~$R$.

Let $I\subset R$ be a non-zero ideal that is also an $O_K$-ideal.
Then $R$ is the inverse image in $O_K$ of the subring $R/I$ of the
quotient $O_K/I$ of~$O_K$). Actually, the diagram:
$$
\begin{CD}
R @>>> R/I \\
@VVV @VVV \\
O_K @>>> O_K/I 
\end{CD}
$$
is both Cartesian and co-Cartesian. For us, the most important is that
every order of $K$ is obtained as follows: take the inverse image in
$O_K$ of a subring of a finite quotient of~$O_K$.

\subsection{Discriminants.}
Recall that $\discr(O_K)$ is the discriminant of the trace form
on~$O_K$.  To be precise: for $M$ a free $\ZZ$-module of finite rank
and $b$ a bilinear form on $M$, let $\discr(M,b)$ be the integer
defined by: let $m$ be a basis of $M$, then $\discr(M,b)$ is the
determinant of the matrix of $b$ relative to~$m$. In more intrinsic
terms, one can use that $b$ induces a bilinear form on the maximal
exterior power of $M$, and use the integer coming from there. The
whole thing does not depend on the basis because changing the basis
changes it by the square of a unit. Over more general rings, and for
projective modules, one obtains an ideal, locally principal, with some
extra structure due to the squares of units that intervene. In fact,
one sees that if the local generators of the ideal are non zero
divisors, then the ideal is, as invertible module, the square of
$\Lambda M$. In our case, we use the bilinear form
$(x,y)\mapsto\trace(xy)$. The separability of $\QQ\to K$ shows that
$\discr(O_K)\neq 0$. Choosing a basis of $O_K$ adapted to $R$ shows
that:
$$
\discr(R)=\discr(O_K)\,|O_K/R|^2. 
$$

\begin{thm}
Let $\zeta_R$ denote the zeta function of the order $R$, i.e., the zeta 
function of $\Spec(R)$ in the usual sense. Then: 
$$
\Res_1(\zeta_R) := \lim_{s\to 1}(s-1)^m\zeta_R(s) = 
\frac{2^{r_1}(2\pi)^{r_2}|\Pic(R)|\,\Reg(R)}
{|\tors(R^*)|\, |\discr(R)|^{1/2}}, 
$$
with $\RR\otimes K\cong\RR^{r_1}\times\CC^{r_2}$ and $\Reg(R)$ the regulator 
of~$R$ (see in the proof for the definition). (Recall that $K$ is the 
product of $m$ number fields.) 
\end{thm}
\begin{proof}
For $O_K$, see for example Lang's ``Algebraic number theory'', 2nd edition, 
VIII,~\S2. In fact, Lang gives the proof when $K$ is a field, but for $O_K$ 
in a product of number fields everything decomposes into products. 
Let us digress a little bit on the regulator. I find that the regulator 
$\Reg'(R)$ should be defined as follows: one considers 
$$
O_K^*\lto (\RR\otimes O_K)^* \overset{\log\|\,\|}{\lto}\RR^{r_1+r_2}, 
$$
and puts: 
$$
\Reg'(R):=\Vol(\RR^{r_1+r_2, +=0}/\text{image of $O_K^*$}), 
$$
with the volume measured with respect to the volume form coming from the 
standard inner product on~$\RR^{r_1+r_2}$, and where $\log\|\,\|$ is taking 
log of absolute value at every factor of $\RR\otimes O_K$, with $\|x\|$ 
being the factor by which the Haar measures change ($|x|$ for a real place, 
$|x|^2$ for a complex one). But this does not give the usual definition, 
as given in Lang. There one omits any one of the infinite places in order to 
get a square matrix of which one takes absolute value of the determinant. 
One easily proves that $\Reg'(R)=2^{-r_2}(r_1+2r_2)(r_1+r_2)^{-1/2}\Reg(R)$, 
which actually makes my definition a bit ugly. 

Anyway, let's proceed. Since we know the theorem for $O_K$, all we have to 
do is to compare our $R$ to~$O_K$. Let $X:=\Spec(O_K)$, $Y:=\Spec(R)$, 
and $N\colon X\to Y$ the morphism induced by the inclusion of $R$ in~$O_K$. 
Then we have a short exact sequence of sheaves on~$X$: 
$$
0\lto \calO_Y^*\lto N_*\calO_{X}^*\lto Q\lto 0, 
$$
with $Q$ a skyscraper sheaf given by $\ol{O_K}^*/\ol{R}^*$ in case $R$ is 
given by the subring $\ol{R}$ of the finite quotient $\ol{O_K}$ of~$O_K$. 
This gives a long exact sequence: 
$$
0\lto R^*\lto O_K^*\lto \ol{O_K}^*/\ol{R}^* \lto\Pic(R)\lto\Pic(O_K)\lto 0. 
$$
Let $A$ be the cokernel of $R^*\to O_K^*$, and $B$ the kernel 
of $\Pic(R)\to\Pic(O_K)$. Then one gets: 
\begin{align*}
&|\Pic(R)| = |B|\,|\Pic(O_K)|, \quad |\ol{O_K}^*/\ol{R}^*| = |A|\,|B|, \quad 
|A| = \frac{|\tors(O_K^*)|}{|\tors(R^*)|}
\left|\frac{O_K^*/\tors}{R^*/\tors}\right| \\
&\Reg(R) = \left|\frac{O_K^*/\tors}{R^*/\tors}\right|\Reg(O_K),\quad 
|\discr(R)|^{1/2} = \frac{|\ol{O_K}|}{|\ol{R}|}|\discr(O_K)|^{1/2}. 
\end{align*}
Putting this all together shows that the right hand side of the equality we 
want to prove changes by the factor 
$|\ol{O_K}^*|\,|\ol{R}^*|^{-1}\,|\ol{R}|\,|\ol{O_K}|^{-1}$ when going from 
$O_K$ to~$R$. 
So all that we have to do now is to show that the left hand side changes 
by the same factor. But then note: 
$$
\frac{|\ol{R}|}{|\ol{R}^*|} = \prod_{\text{$k$ res field of
$\ol{R}$}}\frac{|k|}{|k^*|} = \prod_k \frac{1}{1-|k|^{-1}},
$$
which is clearly the contribution to $\Res_1(\zeta_R)$ of those residue 
fields. 
\end{proof}

\begin{thm}
Let $N>0$. Then there exists a real number $c>0$ such that for every order 
$R$ in a separable $\QQ$-algebra $K$ of degree at most $N$, one has: 
$$
|\Pic(R)|\,\Reg(R) \geq c\,|\discr(R)|^{1/7}. 
$$
\end{thm}
\begin{rem}
As the proof will show, we can actually get $1/6-\eps$ as exponent, 
instead of $1/7$, with a $c$ depending on $\eps$, for every $\eps>0$. 
If one assumes the generalized Riemann hypothesis, then one can get 
$1/2-\eps$ as exponent, for every $\eps>0$, with again $c$ depending 
on~$\eps$. In that case, one uses Siegel's theorem that one finds in 
\cite[Ch.~XIII, \S4]{Lang1}. 
\end{rem}
\begin{proof}
We will first prove this for maximal orders in number fields of bounded 
degree, then for maximal orders in finite separable $\QQ$-algebras of 
bounded degree, and then for arbitrary orders of bounded degree. 

In the case of a maximal order of a number field of bounded degree, we just 
apply two theorems. The first one is the Brauer-Siegel theorem 
(see for example \cite[Ch.~XVI]{Lang1}), that states that: 
\begin{quote} \sl
for $N>0$ and $\eps>0$, there exists $c>0$ such that: 
$$
|\Pic(O_K)|\,\Reg(O_K) \geq c\,|\discr(O_K)|^{1/2-\eps}
$$
for all Galois extensions $K$ of $\QQ$ of degree at most~$N$. 
\end{quote}
The second theorem is one of Stark (\cite[Thm.~1]{Stark1}): 
\begin{quote}\sl 
let $N>0$. There exists $c>0$ such that for all number fields $K$ of 
degree at most $N$ over $\QQ$, one has: 
$
|\Pic(O_K)|\,\Reg(O_K) \geq c\,|\discr(O_K)|^{1/2-1/[K:\QQ]}. 
$
\end{quote}
Together, these two results show:  
\begin{quote}\sl
let $N>0$. There exists $c>0$ such that for every number field $K$ of 
degree at most $N$ over $\QQ$ one has: 
$
|\Pic(O_K)|\,\Reg(O_K) \geq c\,|\discr(O_K)|^{1/6}. 
$
\end{quote}
This settles the case where the $\QQ$-algebra $K$ is a field. The case 
for a maximal order in a finite separable $\QQ$-algebra of degree at most 
$N$ then follows, because everything decomposes into a product of at most 
$N$ factors, for which one has the result already. 

So let now $K$ be a finite separable $\QQ$-algebra of degree at most $N$, and 
let $R$ be an order in it, given by the subring $\ol{R}$ of some finite 
quotient $\ol{O_K}$ of~$O_K$. We have already seen that: 
\begin{align*}
&|\Pic(R)|\,\Reg(R) = \frac{|\ol{O_K}^*|}{|\ol{R}^*|}|\Pic(O_K)|\,\Reg(O_K)
\frac{|\tors(R^*)|}{|\tors(O_K^*)|},\\
&|\discr(R)|=\left(\frac{|\ol{O_K}|}{|\ol{R}|}\right)^2|\discr(O_K)|. 
\end{align*}
We note that the quotient $|\tors(R^*)|\,|\tors(O_K^*)|^{-1}$ and its 
inverse are bounded in terms of $N$ only. Hence the theorem follows from 
the following claim: 
\begin{quote}\sl
for $N>0$ and $\eps>0$ there exists $c>0$ such that for $R$ an order in a 
finite separable $\QQ$-algebra $K$ of degree at most $N$, one has: 
$$
\frac{|\ol{O_K}^*|}{|\ol{R}^*|}\geq 
c\left(\frac{|\ol{O_K}|}{|\ol{R}|}\right)^{1-\eps}, 
$$
where $R$ is the inverse image of the subring $\ol{R}$ of the finite 
quotient $\ol{O_K}$ of~$O_K$. 
\end{quote}
We now prove this claim. Let $n$ denote $|\ol{O_K}/\ol{R}|$. We may and do 
assume that $n>1$. Localizing at the maximal ideals of $\ol{R}$, followed 
by a simple computation, shows that: 
$$
\frac{|\ol{O_K}^*|}{|\ol{R}^*|}\geq n\prod_{p|n}\left(1-\frac{1}{p}\right)^N 
\geq n\left(\frac{1}{5\log(n)}\right)^N. 
$$
Since $\log(n)=n^{o(1)}$, this shows our claim, and hence finishes the 
proof of the theorem. 
\end{proof}

\section{On effective Chebotarev.} 
As usual, let $\Li(x):=\int_2^xdt/\log(t)$. 
If one assumes GRH, then the effective Chebotarev theorem of Lagarias, 
Montgomery and Odlyzko, stated as in~\cite[Thm.~4]{Serre1} and the 
second remark following that theorem, says: 
\begin{quote}\sl
for $M$ a finite Galois extension of $\QQ$, let $n_M$ denote its degree, 
$d_M$ its absolute discriminant $|\discr(O_M)|$, and for $x$ in $\RR$, let 
$\pi_{M,1}(x)$ be the number of primes $p\leq x$ that are unramified in $M$ 
and such that the Frobenius conjugacy class $\Frob_p$ contains just the 
identity element of~$\Gal(M/\QQ)$. Then one has, for all sufficiently large 
$x$ and all finite Galois extensions $M$ of~$\QQ$: 
$$
\left|\pi_{M,1}(x)-\frac{1}{n_M}\Li(x)\right| \leq 
\frac{1}{3n_M}x^{1/2}\left(\log(d_M) + n_M\log(x)\right). 
$$
\end{quote}
This result shows that for all $x$ sufficiently large, and all finite Galois 
extensions $M$ of $\QQ$, one has: 
$$
\pi_{M,1}(x) \geq \frac{x}{n_M\log(x)}
\left(\Li(x)\frac{\log(x)}{x} - \frac{\log(x)}{3x^{1/2}} 
\left(\log(d_M)+n_M\log(x)\right)\right).
$$
If $x$ tends to infinity, $\Li(x)\log(x)/x$ tends to 1 and 
$\log(x)^2/x^{1/2}$ tends to~0. Some computation (that we will do below) 
shows that if $x$ is sufficiently big (i.e., bigger than some absolute 
constant), and bigger than $2(\log(d_M)^2(\log(\log(d_M)))^2$, then 
$\log(x)\log(d_M)/3x^{1/2}<1/2$, and hence: 
$$
\pi_{M,1}(x)\geq \frac{x}{3n_M\log(x)}. 
$$
Here is the computation that I promised. Put $a:=\log(d_M)/3$. We want to 
find a lower bound for $x$ that implies that $a\log(x)/\sqrt{x}<1/2$. 
We put $x:=y^2$ (with $y>0$, of course). Then what we want is a lower bound 
for $y$ such that $b\log(y)<y$, with $b=4a$. we put $y=zb$. Then what we 
want is a lower bound for $z$ such that $z-\log(z)>c$, with $c=\log(b)$. 
Now write $z=(1+u)c$. Then what we want is: $uc-\log(1+u)-\log(c)>0$. 
Since $\log(1+u)\leq u$, it suffices that $uc-u-\log(c)>0$, i.e., 
that $u>\log(c)/(c-1)$ (by the way, since we are willing to let $x$ be 
sufficiently large, we may take care of small $d_M$ by that, and suppose 
that $c$ is sufficiently large). For $c$ sufficiently large, for any $\eps>0$, 
$u>\eps$ is good enough. Translating this back to $x$ and $\log(d_M)$, 
one gets that $x>(1+\eps)^2(4a\log(4a))^2$ is good. Then one uses that 
$2>16/9$. 

\section{Real approximation.}
It is known that for $G$ an affine algebraic group over $\QQ$ one has 
$G(\QQ)$ dense in~$G(\RR)$. This is what Deligne calls real approximation. 
To prove it, one reduces to the case of tori. But even that case is not so 
trivial to me. Of course, tori are unirational (they are images of tori that 
are products of multiplicative groups of number fields), but that is not 
enough: that only gives that the rational points of $G$ are dense in the 
connected components of $G(\RR)$ that do contain a rational point. Anyway, 
for a detailed proof I would refer to the book \cite{Platonov-Rapinchuk} 
of Platonov  and Rapinchuk. 

\end{document}